\documentclass{article}

\usepackage{amssymb,latexsym,amsmath}

\usepackage{graphicx}

\begin{document}

\newcommand{\bfi}{\bfseries\itshape}

\makeatletter

\@addtoreset{figure}{section}

\def\thefigure{\thesection.\@arabic\c@figure}

\def\fps@figure{h, t}

\@addtoreset{table}{bsection}

\def\thetable{\thesection.\@arabic\c@table}

\def\fps@table{h, t}

\@addtoreset{equation}{section}

\def\theequation{\thesubsection.\arabic{equation}}

\makeatother

\newtheorem{thm}{Theorem}[section]

\newtheorem{prop}[thm]{Proposition}

\newtheorem{lem}[thm]{Lemma}

\newtheorem{cor}[thm]{Corollary}

\newtheorem{dfn}[thm]{Definition}

\newtheorem{rmk}[thm]{Remark}

\newtheorem{exempl}{Example}[section]

\newenvironment{exemplu}{\begin{exempl}  \em}{\hfill $\square$

\end{exempl}}

\newcommand{\comment}[1]{\par\noindent{\raggedright\texttt{#1}

\par\marginpar{\textsc{Comment}}}}

\newcommand{\todo}[1]{\vspace{5 mm}\par \noindent \marginpar{\textsc{ToDo}}\framebox{\begin{minipage}[c]{0.95 \textwidth}

\tt #1 \end{minipage}}\vspace{5 mm}\par}

\newcommand{\ea}{\mbox{{\bf a}}}

\newcommand{\eu}{\mbox{{\bf u}}}

\newcommand{\ueu}{\underline{\eu}}

\newcommand{\ueo}{\overline{u}}

\newcommand{\oeu}{\overline{\eu}}

\newcommand{\ew}{\mbox{{\bf w}}}

\newcommand{\ef}{\mbox{{\bf f}}} 

\newcommand{\eF}{\mbox{{\bf F}}}

\newcommand{\eC}{\mbox{{\bf C}}}

\newcommand{\en}{\mbox{{\bf n}}}

\newcommand{\eT}{\mbox{{\bf T}}}

\newcommand{\eL}{\mbox{{\bf L}}}

\newcommand{\eR}{\mbox{{\bf R}}}

\newcommand{\eV}{\mbox{{\bf V}}}

\newcommand{\eU}{\mbox{{\bf U}}}

\newcommand{\ev}{\mbox{{\bf v}}}

\newcommand{\eve}{\mbox{{\bf e}}}

\newcommand{\uev}{\underline{\ev}}

\newcommand{\eY}{\mbox{{\bf Y}}}

\newcommand{\eK}{\mbox{{\bf K}}}

\newcommand{\eP}{\mbox{{\bf P}}}

\newcommand{\eS}{\mbox{{\bf S}}}

\newcommand{\eJ}{\mbox{{\bf J}}}

\newcommand{\eB}{\mbox{{\bf B}}}

\newcommand{\eH}{\mbox{{\bf H}}}

\newcommand{\leb}{\mathcal{ L}^{n}}

\newcommand{\eI}{\mathcal{ I}}

\newcommand{\eE}{\mathcal{ E}}

\newcommand{\hen}{\mathcal{H}^{n-1}}

\newcommand{\eBV}{\mbox{{\bf BV}}}

\newcommand{\eA}{\mbox{{\bf A}}}

\newcommand{\eSBV}{\mbox{{\bf SBV}}}

\newcommand{\eBD}{\mbox{{\bf BD}}}

\newcommand{\eSBD}{\mbox{{\bf SBD}}}

\newcommand{\ecs}{\mbox{{\bf X}}}

\newcommand{\eg}{\mbox{{\bf g}}}

\newcommand{\paromega}{\partial \Omega}

\newcommand{\gau}{\Gamma_{u}}

\newcommand{\gaf}{\Gamma_{f}}

\newcommand{\sig}{{\bf \sigma}}

\newcommand{\gac}{\Gamma_{\mbox{{\bf c}}}}

\newcommand{\deu}{\dot{\eu}}

\newcommand{\dueu}{\underline{\deu}}

\newcommand{\dev}{\dot{\ev}}

\newcommand{\duev}{\underline{\dev}}

\newcommand{\weak}{\stackrel{w}{\approx}}

\newcommand{\mild}{\stackrel{m}{\approx}}

\newcommand{\strong}{\stackrel{s}{\approx}}

\newcommand{\weakdown}{\rightharpoondown}

\newcommand{\opg}{\stackrel{\mathfrak{g}}{\cdot}}

\newcommand{\opunu}{\stackrel{1}{\cdot}}
\newcommand{\opdoi}{\stackrel{2}{\cdot}}

\newcommand{\opn}{\stackrel{\mathfrak{n}}{\cdot}}
\newcommand{\opx}{\stackrel{x}{\cdot}}

\newcommand{\tr}{\ \mbox{tr}}

\newcommand{\Ad}{\ \mbox{Ad}}

\newcommand{\ad}{\ \mbox{ad}}

\title{Infinitesimal affine geometry of metric spaces endowed with a dilatation structure}

\author{Marius Buliga \\
\\
Institute of Mathematics, Romanian Academy \\
P.O. BOX 1-764, RO 014700\\
Bucure\c sti, Romania\\
{\footnotesize Marius.Buliga@imar.ro}}

\date{This version:  31.03.2008}

\maketitle

\begin{abstract}
We study algebraic and geometric properties of metric spaces endowed with 
dilatation structures, which are emergent during the passage through smaller and smaller scales. In the limit 
we obtain  a generalization of metric affine geometry,  endowed with  a  noncommutative vector addition operation and with 
a modified version of ratio of three collinear points. This is the geometry of normed affine group spaces, a category which 
contains the ones  of  homogeneous groups, Carnot groups or contractible groups. In  this category group operations are 
not fundamental, but derived objects, and the generalization of affine geometry is not based on incidence relations. 

\end{abstract}


\paragraph{MSC classes:} 	20F65; 20F19;  22A10

\newpage

\tableofcontents

\section{Introduction}

The point of view that dilatations can be taken as fundamental objects which 
 induce a differential calculus is relatively well known. The idea is simple: 
in a vector space $\mathbb{V}$ define the dilatation based at $x$ and of coefficient 
$\varepsilon > 0$ as the function which associates to $y$ the value 
$$\delta^{x}_{\varepsilon} y \ = \ x + \varepsilon (y - x) \quad . $$
Then for  a function $f: \mathbb{V} \rightarrow \mathbb{W}$  
between vector spaces $\mathbb{V}$ and $\mathbb{W}$ we have: 
$$ \left( \delta^{f(x)}_{\varepsilon^{-1}}  f  \delta^{x}_{\varepsilon} \right) (u) \ = \ f(x) + \frac{1}{\varepsilon} \left[ f(x + \varepsilon (u-x)) - f(x) \right]  \quad , $$
thus the directional derivative of $f$ at $x$, along $u-x$ appears as: 
$$f(x) \, + \, D \, f(x) (u-x) \ = \ \lim_{\varepsilon \rightarrow 0}   \left( \delta^{f(x)}_{\varepsilon^{-1}}  f  \delta^{x}_{\varepsilon} \right) (u) \quad . $$
Until recently there was not much interest into the generalization of such a differential calculus, based on other 
dilatations than the usual ones, probably 
because  nobody knew any fundamentally different  example. 

This changed gradually due to different lines of research, like the study of hypoelliptic operators 
H\"ormander \cite{hormander},  harmonic analysis on homogeneous groups 
Folland, Stein \cite{fostein},  probability theory on groups Hazod \cite{hazod}, Siebert \cite{siebert}, studies in geometric 
analysis in metric spaces  in relation with
sub-riemannian geometry   Bella\"{\i}che   \cite{bell}, groups with 
polynomial growth   Gromov \cite{gromovgr}, or Margulis type rigidity  
results   Pansu \cite{pansu}. 

Another line of research concerns the differential calculus 
over general base fields and rings,   Bertram, Gl\"{o}ckner and Neeb \cite{bertram2}. As the authors explain, 
it is possible to construct such a differential calculus without using the specific properties of the base field 
(or ring). In their approach it is not made a distinction between real and ultrametric differential calculus (and 
even not between finite dimensional and infinite dimensional differential calculus). They point out that differential 
calculus (integral calculus not included) seems to be a part of analysis which is completely general, based only on 
elementary results in linear algebra and topology. 

The differential calculus proposed by Bertram, Gl\"{o}ckner and Neeb is  a generalization of ``classical'' 
calculus in topological vector spaces over general base fields, and even over rings. The operation of vector addition 
is therefore abelian, modifications being made in relation with the multiplication by scalars.

A different idea, emergent in the studies concerning geometric analysis in metric spaces, is to establish a differential 
calculus in homogeneous groups, in particular in Carnot groups. These are noncommutative versions of topological vector 
spaces, in the sense that the operation of addition (of ``vectors'') is replaced by a noncommutative group operation 
and there is a replacement of multiplication by scalars in a general base field with a multiplicative action of $(0, + \infty)$ by 
group automorphisms. 

In fact this is only a part of the nonsmooth calculus encountered in geometric analysis on metric spaces. For a survey see the paper by 
Heinonen \cite{heinonen}. The objects of interest in nonsmooth calculus as described by Heinonen are spaces of homogeneous type, or 
metric measured spaces where a generalization of Poincar\'e inequality  is true. In such spaces the differential calculus goes a long way: 
Sobolev spaces, differentiation theorems, Hardy spaces. It is noticeable that in such a general situation we don't have  enough structure 
to define differentials, but only various constructions corresponding to the norm of a differential of a function. Nevertheless 
see the remarkable result of Cheeger \cite{cheeger}, who proves that to  a metric measure space satisfying a Poincar\'e inequality we can 
associate an $L^{\infty}$ cotangent bundle with finite dimensional fibers. Other important works which might also be relevant in relation to this 
paper are David, Semmes \cite{davsemmes}, where spaces with arbitrary small neighbouhoods containing similar images of the whole space are studied, 
or David, Semmes \cite{davsemmes2}, where they study rectifiability properties of subsets of $\mathbb{R}^{n}$  with arbitrary small neighbourhoods containing ``big pieces of bi-Lipschitz images'' of the whole subset.

A particular case of a space of homogeneous type where more can be said is a normed homogeneous group, definition \ref{defnormedhom}.  
According to \cite{fostein} p. 5, a homogeneous group is a connected and simply connected Lie group whose Lie algebra is endowed with a 
family of dilatations $\displaystyle \left\{ \delta_{\varepsilon} \mbox{ : } \varepsilon \in (0, + \infty) \right\}$, which are algebra automorphisms, 
simultaneously diagonalizable. As in this case the exponential of the group is a bijective mapping, we may transform dilatations of the algebra 
into dilatations of the group, therefore homogeneous groups are conical groups. Also, they can be described as nilpotent Lie groups positively 
graded. 

Carnot groups are homogeneous groups which are stratified, meaning that the first nontrivial element of the graduation generates the whole group (or algebra). 
The interest into such groups come from various sources, related 
mainly to  the study of hypo-elliptic operators   H\"ormander \cite{hormander}, and to extensions 
of  harmonic analysis   Folland, Stein \cite{fostein}. 

Pansu introduced the first really new example of such a  differential calculus based on other than 
usual dilatations: the ones which are associated to a Carnot group. He proved in \cite{pansu} the potential of what is 
now called Pansu derivative, by providing an alternative proof of a Margulis rigidity type result, as a corollary  
of the Rademacher theorem for Lipschitz functions on Carnot groups. 
Rademacher theorem, stating that a Lipschitz function is derivable 
almost everywhere, is a mathematical crossroad, because there meet measure theory, differential calculus and metric geometry. 
In \cite{pansu} Pansu proves a generalization of this theorem for his new derivative. 

The challenge to extend Pansu results to general regular sub-riemannian manifolds, taken by Margulis, Mostow \cite{marmos1} \cite{marmos2}, Vodopyanov \cite{vodopis} and others, is difficult because on such general metric space there is no natural underlying algebraic structure, as in the case of 
Carnot groups, where we have the group operation as a non commutative replacement of the operation of addition in vector spaces. 

On a regular sub-riemannian manifold we have  
to construct simultaneously several objects: tangent spaces to a point in the sub-riemannian space, an operation of addition of ``vectors'' in the 
tangent space, and a derivative of the type considered by Pansu. Dedicated to the first two objects is a string of papers, either directly related 
to the subject, as Bella\"{\i}che \cite{bell}, or growing on techniques which appeared in the paper dedicated to groups of polynomial 
growth of Gromov \cite{gromovgr}, continuing in the big paper Gromov \cite{gromovsr}. 

In these papers dedicated to sub-riemannian geometry the lack of a underlying algebraic structure was supplanted by using techniques of differential 
geometry. At a closer look, this means that in order to construct the fundamentals of a non standard differential calculus, the authors used the 
classical one. This seems to me  comparable to efforts to study hyperbolic geometry on models, like the
Poincar\'e disk, instead of intrinsically 
explore the said geometry. 

Dilatation structures on  metric spaces,  introduced in \cite{buligadil1},  
describe the approximate self-similarity properties of a metric space. A 
dilatation structure is a notion related, but more general, to groups  and  
 differential structures.  

The basic objects of a dilatation structure are dilatations (or contractions). 
The axioms of a dilatation structure set the rules of interaction between
different dilatations.

The point of view of dilatation structures is that dilatations are really fundamental objects, not only for 
defining a notion of derivative, but as well for all algebraic structures that we may need.

This viewpoint is justified by the following results obtained in \cite{buligadil1}, explained 
in a condensed and improved presentation, in the first part of this paper.  
A metric space $(X,d)$ which admits a strong dilatation structure (definition 
\ref{defweakstrong}) has a metric tangent space at any point $x \in X$ (theorem 
\ref{thcone}), and any such metric tangent space has an algebraic structure 
of a conical group (theorem \ref{tgene}). 

Conical groups are generalizations of homogeneous Lie groups, but also of p-adic nilpotent groups, or of general contractible groups. 
A conical group is a locally compact group endowed with a family of dilatations 
$\displaystyle \left\{ \delta_{\varepsilon} \mbox{ : } \varepsilon \in \Gamma \right\}$. 
 Here  $\Gamma$ is a locally compact abelian group with an associated  morphism $\nu: \Gamma \rightarrow (0, + \infty)$ which 
distinguishes an end of $\Gamma$, namely the filter generated by the pre-images $\displaystyle \nu^{-1}(0, r)$, $r > 0$. This end, is denoted 
by $0$ and $\varepsilon \in \Gamma \, \rightarrow 0$ means $\nu(\varepsilon) \rightarrow 0$ in $(0, + \infty)$. 
 Any contractible group is a conical group and to any conical group 
we can associate a family of contractible groups.

 The structure of contractible groups is known in some detail, due to Siebert
\cite{siebert}, Wang \cite{wang}, Gl\"{o}ckner and Willis \cite{glockwill}, 
Gl\"{o}ckner \cite{glockner} and references therein.

By a classical result of Siebert \cite{siebert}
proposition 5.4, we can characterize the algebraic structure of the metric
tangent spaces associated to dilatation structures of a certain kind: they are  homogeneous groups  
(corollary \ref{cortang}). The corollary \ref{cortang} thus represents a
generalization  of difficult results in sub-riemannian geometry concerning the structure of the 
metric tangent space at a point of a regular sub-riemannian manifold. This line of research is pursued 
further in the paper \cite{buligaradon}. 
 
Morphisms of dilatation structures generalize the notion of affine transformation. 
A dilatation structure on a metric space induce a family of dilatation structures 
on the same space, at different scales. We explain that canonical morphisms between these 
induced dilatation structures lead us to a kind of emergent affinity on smaller and 
smaller scale. 

Finally we  characterize  contractible groups in terms of
dilatation structures. To a normed  contractible group we 
can naturally associate a linear dilatation structure (proposition
\ref{pexlin}). Conversely, by theorem  \ref{tdilatlin} any linear and strong 
dilatation structure comes  from the  dilatation structure of a 
normed contractible group.

We are thus led to the introduction of a noncommutative affine geometry, in the spirit of 
Bertram ``affine algebra'', which is commutative according to our point of view. In such a geometry 
incidence relations are no longer relevant, being  replaced by algebraic axioms concerning dilatations. 
We define a version of the ratio of three collinear points (replaced by a ``ratio function'' which associates 
to a pair of points and two positive numbers the third point) and we prove that it is the basic invariant of 
this geometry. Moreover, it turns out that this is the geometry of normed affine group spaces, a notion which 
is to conical groups as a normed affine space is to a normed topological vector space (theorem \ref{taffine}).

\section{Affine structure in terms of dilatations}

\subsection{Affine algebra}

Bertram \cite{bertram1} Theorem 1.1 (here theorem \ref{tbertram}) and paragraph 5.2, proposes the following algebraic description of affine geometry and of affine metric geometry over a field $\mathbb{K}$ of characteristic different from $2$, which is not based on incidence notions, but on algebraic relations concerning ``product maps''. He then pursues to the development of generalized projective geometries and their relations 
to Jordan algebras. 
For our purposes,  
we changed the name of ``product maps'' (see the theorem below) from ``$\displaystyle \pi$'' to ``$\displaystyle \delta$'', more precisely: 
$$\pi_{r}(x,y) \ = \ \delta_{r}^{x} y $$
Further, in theorem \ref{tbertram} and definition \ref{dbertram} is explained this point of view. 

\begin{thm}
 The category of affine spaces over a field $\mathbb{K}$ of characteristic different from $2$  is equivalent with 
the category of sets $\mathbb{V}$ equipped with a family $\displaystyle \delta_{r}$, $r \in \mathbb{K}$, of ``product 
maps''
$$\delta_{r}: \mathbb{V} \times \mathbb{V} \rightarrow \mathbb{V} \, , \ \ (x,y) \mapsto \delta_{r}^{x} y $$
satisfying the following properties (Af1) - (Af4): 
\begin{enumerate}
 \item[(Af1)] The map $\displaystyle r \mapsto \delta_{r}^{x}$ is a homomorphism of the unit group $\displaystyle \mathbb{K}^{\times}$ into the group 
of bijections of $\mathbb{V}$ fixing $x$, that is 
$$\delta_{1}^{x} y = y \ , \ \delta_{r}^{x} \,  \delta_{s}^{x} y = \delta_{rs}^{x} y \ , \ \delta_{r}^{x} x = x $$
\item[(Af2)] For all $r \in \mathbb{K}$ and $x \in \mathbb{V}$ the map
$\displaystyle \delta^{x}_{r}$ is an endomorphism of $\displaystyle \delta_{s}$, 
$s \in \mathbb{K}$: 
$$\delta_{r}^{x} \,  \delta_{s}^{y} z = \delta_{s}^{\delta_{r}^{x} y}  \delta_{r}^{x} z$$
\item[(Af3)] The ``barycentric condition'': 
 $\displaystyle \delta_{r}^{x} \, y \, = \delta_{1-r}^{y} \, x$
\item[(Af4)] The group generated by the $\displaystyle \delta_{r}^{x}  \, \delta_{r^{-1}}^{y}$ ($\displaystyle r \in \mathbb{K}^{\times}$, $x,y \in \mathbb{V}$) is abelian, that 
is for all $\displaystyle r, s \in \mathbb{K}^{\times}$, $x, y, u, v \in \mathbb{V}$ 
$$ \delta_{r}^{x}  \, \delta_{r^{-1}}^{y} \,   \delta_{s}^{u}  \, \delta_{s^{-1}}^{v} \ = \ \delta_{s}^{u}  \, \delta_{s^{-1}}^{v} \, \delta_{r}^{x}  \, \delta_{r^{-1}}^{y}$$
\end{enumerate}
More precisely, in every affine space over $\mathbb{K}$, the maps 
\begin{equation}
 \delta_{r}^{x} y = (1-r) x + r y
\end{equation}
with $r \in \mathbb{K}$, satisfy (Af1) - (Af4). Conversely, if product maps with the properties (Af1) - (Af4) are given and $x\in \mathbb{V}$ is 
an arbitrary point then 
$$ u \, +_{x} \, v : = \delta_{2}^{x} \,  \delta_{\frac{1}{2}}^{u} v \quad , \quad r \, u \, : = \delta_{r}^{x} u$$
defines on $\mathbb{V}$ the structure of a vector space over $\mathbb{K}$ with zero vector $x$, and this construction is inverse to the 
preceding one. Affine maps $g: \mathbb{V} \rightarrow \mathbb{V}'$ in the usual sense are precisely the homomorphisms of product maps, that is 
maps $g: \mathbb{V} \rightarrow \mathbb{V}'$ such that $\displaystyle g \, \pi_{r}(x,y) \, = \, \pi_{r}'( g x , g y)$ for all 
$x,y \in \mathbb{V}$, $r \in \mathbb{K}$. 
\label{tbertram}
\end{thm}

We shall use the name ``real normed affine space'' in the following sense. 

\begin{dfn}
 A real normed affine space is an affine space $\mathbb{V}$ over $\mathbb{R}$ together with a distance function $d: \mathbb{V} \times \mathbb{V} \rightarrow \mathbb{K}$ such that: 
\begin{enumerate}
 \item[(Af5)] for all $x \in \mathbb{V}$ $\displaystyle \| \cdot \|_{x}: = d(x, \cdot) : \mathbb{V} \rightarrow \mathbb{K}$ is a  norm on the vector space 
$(\mathbb{V}, x)$ with zero vector $x$. 
\item[(Af6)] the distance $d$   is translation invariant: for any $x, y, u, v \in \mathbb{V}$ we have: 
$$d(x +_{u} v, y +_{u} v) = d(x,y)$$ 
\end{enumerate}
\label{dbertram}
\end{dfn}

We remark that the field of product maps $\displaystyle \delta^{x}_{r}$ (together with the distance function $d$ for the metric 
case) is the central object in the construction of affine geometry over a general field.

\subsection{Focus on dilatations}

There is another, but related, way of generalizing the affine geometry, which is the one of dilatation structures 
\cite{buligadil1}. In this approach product maps of Bertram are replaced by ``dilatations''. 

For this  we have to   replace the field $\mathbb{K}$ by a commutative group $\Gamma$ (instead of the 
multiplicative group $\displaystyle \mathbb{K}^{\times}$) endowed with a ``valuation map'' $\nu: \Gamma \rightarrow (0,+\infty)$, 
which is a group morphism. We  write $\varepsilon \rightarrow 0$, $\varepsilon \in \Gamma$, for $\nu(\varepsilon) \rightarrow 0$ 
in $(0, +\infty)$. We  keep axioms like (Af1), (Af2) (from Theorem \ref{tbertram}), but we
 modify  (Af5) (from Definition \ref{dbertram}). There will be one more axiom concerning the relations 
between the distance and dilatations. This is explained in theorem \ref{taffine}. 

The conditions appearing in theorem \ref{taffine} are a particular case of the system of axioms of dilatation structures, 
introduced in \cite{buligadil1}. Dilatation structures are also a generalization of homogeneous groups, definition \ref{defnormedhom}, in fact 
we arrived to dilatation structures after an effort to find a common algebraic and analytical ground for homogeneous groups and 
sub-riemannian manifolds.

The axioms of a dilatation structure are partly algebraic and partly of an analytical nature (by using uniform limits). Metric spaces 
endowed with dilatation structures have beautiful properties. The most important is that for any point in such a space there is a tangent 
space (in the metric sense) realized as a ``normed conical group''. Any normed conical group has an associated dilatation structure 
which is   ``linear'' in the sense that it satisfies (Af2). However, conical groups form a family much larger than affine spaces (in the 
usual sense, over $\mathbb{R}$ or $\mathbb{C}$). Building blocks of conical groups are homogeneous groups (graded Lie groups) or 
p-adic versions of them. By renouncing to (Af3) and (Af4) we thus allow noncommutativity of the ``vector addition'' operation.

Let us explain how we can recover the usual affine geometry from the viewpoint of dilatation structures. For simplicity we take here 
$\Gamma = (0,+\infty)$ and $\mathbb{V}$ is a real, finite dimensional vector space.

Here is  the definition of a normed homogeneous group. See section \ref{carnotgroups} for more details on the particular case of stratified homogeneous 
groups. 

\begin{dfn}
 A normed  homogeneous group is a connected and simply connected Lie group whose Lie algebra is endowed with a 
family of dilatations $\displaystyle \left\{ \delta_{\varepsilon} \mbox{ : } \varepsilon \in (0, + \infty) \right\}$, which are algebra automorphisms, 
simultaneously diagonalizable, together with a homogeneous norm. 

Since the Lie group exponential is a bijection we shall identify the Lie algebra with the Lie group, thus a normed homogeneous group  is a group 
operation on a finite dimensional vector space $\mathbb{V}$. The operation will be denoted multiplicatively, with $0$ as neutral element, as in 
Folland, Stein \cite{fostein}. We thus have    a linear action  $\delta: (0,+\infty) \rightarrow Lin(\mathbb{V}, \mathbb{V})$ on $\mathbb{V}$, and a homogeneous norm $\| \cdot \|: \mathbb{V} \rightarrow [0, + \infty)$,  such that: 
\begin{enumerate}
 \item[(a)] for any $\varepsilon \in (0, + \infty)$ the transformation $\displaystyle \delta_{\varepsilon}$ is an automorphism of the group operation: 
for any $x, y \in \mathbb{V}$ we have $\displaystyle \delta_{\varepsilon}( x \cdot y ) \ = \ \delta_{\varepsilon} x \cdot \delta_{\varepsilon} y$ 
\item[(b)] the family $\displaystyle \left\{ \delta_{\varepsilon} \mbox{ : } \varepsilon \in (0, + \infty) \right\}$ is simultaneously diagonalizable: there 
is a finite direct sum decomposition of the vector space $\mathbb{V}$ 
$$\mathbb{V} \ = \ V_{1} + ... + V_{m}$$ 
such that for any $\varepsilon \in (0, + \infty)$ we have: 
$$ x = \sum_{i=1}^{m} x_{i} \, \in \mathbb{V} \,m \mapsto \, \delta_{\varepsilon} x \, = \, \sum_{i=1}^{m}
\varepsilon^{i}  x_{i} \quad .$$
\item[(c)]  the homogeneous norm has the properties: 
\begin{enumerate}
\item[(c1)] $\| x \| = 0$ if and only if $x = 0$, 
\item[(c2)] $\| x \cdot y \| \leq \|x \| + \|y \|$ for any $x, y \in \mathbb{V}$, 
\item[(c3)] for any $x \in \mathbb{V}$ and $\varepsilon > 0$ we have $\displaystyle \| \delta_{\varepsilon} x \| \ = \ \varepsilon \, \|x\|$
\end{enumerate}
\end{enumerate}
\label{defnormedhom}
\end{dfn}

\begin{dfn}
To a normed homogeneous group $\displaystyle (\mathbb{V}, \delta, \cdot , \| \cdot \|)$ we associate a normed affine group space $(\mathbb{V}, +_{\cdot} , \delta^{\cdot}_{\cdot}, d)$. Here  we  use the sign ``$+$'' for an operation which was denoted multiplicatively, 
for compatibility with the previous approach of Bertram , see theorem \ref{tbertram}. The normed affine group space $(\mathbb{V}, +_{\cdot} , \delta^{\cdot}_{\cdot}, d)$ is described by the following points: 
\begin{enumerate}
\item[-] for any $ u  \in \mathbb{V}$ the function $\displaystyle +_{u}: \mathbb{V} \times \mathbb{V} \rightarrow \mathbb{V}$, $\displaystyle x+_{u}v = x \cdot u^{-1} \cdot v$ 
is the left translation of the group operation $\cdot$ with the zero element  $u$. In particular we have 
$\displaystyle x+_{0}y = x \cdot y$.  
\item[-] for any $x, y \in \mathbb{V}$ and $\varepsilon \in (0,+\infty)$ we define  
$$\delta_{\varepsilon}^{x} y \ = \ x \cdot \delta_{\varepsilon} (x^{-1}\cdot y)$$
and remark that the definition is invariant with the choice of the base point for the operation in the sense: for any $u \in \mathbb{V}$ we have: 
$$\delta_{\varepsilon}^{x} y \ = \ x +_{u} \delta_{\varepsilon}^{u} \left( \, inv^{u}(x) +_{u} y \right)$$
where $\displaystyle \, inv^{u}(x)$ is the inverse of $x$ with respect to the operation $\displaystyle +_{u}$, (by computation we get 
$\displaystyle \, inv^{u} (x) = u \cdot x^{-1} \cdot u$), 
\item[-] the distance $d$ is defined as: for any $x, y \in \mathbb{V}$ we have $\displaystyle d(x,y) = \| x^{-1} \cdot y \|$. As previously, remark that 
the definition does not depend on the choice of the base point for the operation, that is: for any $u \in \mathbb{V}$ we have 
$$d(x,y) \ = \ \| inv^{u}(x) +_{u} y \|_{u} \, , \ \| x \|_{u} : = \| u^{-1} \cdot x \|$$
Equally, this is a consequence of the invariance of the norm with respect to left translations (by any group operation $\displaystyle +_{u}$, $u \in \mathbb{V}$). 
\end{enumerate}
\label{defnormedhom1}
\end{dfn}

\begin{thm}
 The category of normed affine group spaces   is equivalent with  
the category of  locally compact metric spaces  $(X,d)$ equipped with a family $\displaystyle \delta_{\varepsilon}$, 
$\varepsilon \in (0,+\infty)$, of dilatations
$$\delta_{\varepsilon}: X \times X \rightarrow X \, , \ \ (x,y) \mapsto \delta_{\varepsilon}^{x} y $$
satisfying the following properties: 
\begin{enumerate}
 \item[(Af1')] The map $\displaystyle \varepsilon \mapsto \delta{\varepsilon}^{x}$ is a homomorphism of the multiplicative group $\displaystyle (0, + \infty)$ into the group 
of continuous, with continuous inverse  functions of $X$ fixing $x$, that is 
$$\delta_{1}^{x} y = y \ , \ \delta_{r}^{x} \,  \delta_{s}^{x} y = \delta_{rs}^{x} y \ , \ \delta_{r}^{x} x = x $$
\item[(A2)] the function $\displaystyle \delta : (0, +\infty) \times X \times X \rightarrow  X$ defined by 
$\displaystyle \delta (\varepsilon,  x, y)  = \delta^{x}_{\varepsilon} y$ is continuous. Moreover, it can be 
continuously extended to $\displaystyle [0,+\infty) \times X \times X$ by $\delta (0,x, y) = x$  and the limit  
$$\lim_{\varepsilon\rightarrow 0} \delta_{\varepsilon}^{x} y \, = \, x  $$
is uniform with respect to $x,y$ in compact set.
\item[(A3')] for any $x \in X$ and for any $u, v \in X$, $\varepsilon \in (0, + \infty)$ we have 
$$\frac{1}{\varepsilon} \, d \left( \delta^{x}_{\varepsilon} u , \delta^{x}_{\varepsilon} v \right) \ = \ d(u,v)$$

\item[(A4)] for any $x, u, v \in X$, $\varepsilon \in (0, + \infty)$ let us define 
$$\Delta^{x}_{\varepsilon}(u,v) = \delta_{\varepsilon^{-1}}^{\delta^{x}_{\varepsilon} u} \delta^{x}_{\varepsilon} v . $$
Then we have the limit 
$$\lim_{\varepsilon \rightarrow 0}  \Delta^{x}_{\varepsilon}(u,v) =  \Delta^{x}(u, v)  $$
uniformly with respect to $x, u, v$ in compact set. 
\item[(Af2')] For all $\varepsilon \in (0,+\infty)$ and $x \in X$ the map $\displaystyle \delta_{\varepsilon}^{x}$ is an endomorphism of $\displaystyle \delta_{s}$, 
$s \in (0,+\infty)$: 
$$\delta_{r}^{x} \,  \delta_{s}^{y} z = \delta_{s}^{\delta_{r}^{x} y}  \delta_{r}^{x} z$$
\end{enumerate}
More precisely, in every normed affine group space , the maps $\displaystyle \delta_{\varepsilon}^{x}$
 and distance $d$ satisfy (Af1'), (A2), (A3'), (A4), (Af2'). Conversely, if dilatations $\displaystyle \delta_{\varepsilon}^{x}$ and distance  $d$ 
are given, such that they satisfy the collection (Af1'), (A2), (A3'), (A4), (Af2'), for an arbitrary point  $x\in \mathbb{V}$  the following expression  
$$ \Sigma^{x}(u,v) : =  \lim_{\varepsilon \rightarrow 0} \delta_{\varepsilon^{-1}}^{x} \,  \delta_{\varepsilon}^{\delta_{\varepsilon}^{x} u } v$$
together with $\displaystyle \delta_{\varepsilon}^{x}$ and distance  $d$ 
defines on $\mathbb{V}$ the structure of a normed affine group space, and this construction is inverse to the 
preceding one. The arrows of this category are bilipschitz invertible  homomorphisms of dilatations, that is 
maps $g: \mathbb{V} \rightarrow \hat{\mathbb{V}}$ such that $\displaystyle g \, \delta_{\varepsilon}^{x} y  \, = \, \hat{\delta}_{r}^{g x} g y$ for all 
$x,y \in \mathbb{V}$, $\varepsilon \in (0,+\infty)$.

Moreover, the category of real normed affine spaces is a subcategory of the previous one, namely the category of  locally compact metric spaces  $(X,d)$ equipped with a family $\displaystyle \delta_{\varepsilon}$, 
$\varepsilon \in (0,+\infty)$, of dilatations satisfying (Af1'), (A2), (A3'), (A4), (Af2') and
\begin{enumerate}
 \item[(Af3)] the ``barycentric condition'': for all $\varepsilon \in (0,1)$ 
 $\displaystyle \delta_{\varepsilon}^{x} \, y \, = \delta_{1-\varepsilon}^{y} \, x$
\end{enumerate}
The arrows of this category are exactly the affine, invertible maps. 
 \label{taffine}
\end{thm}

\paragraph{Proof.} 
Here we shall prove the easy implication, namely  why the conditions (Af1'), (A2), (A3'), (A4), 
(Af2') and (Af3) are satisfied in a real normed affine space.

 For the real normed affine space  space $\mathbb{V}$ let us fix for simplicity a point $0 \in \mathbb{V}$ and work with the vector space $\mathbb{V}$ with 
zero vector $0$. Since a real normed affine space is a particular example of a homogeneous group, definition \ref{defnormedhom} and observations inside apply. 
The dilatation based 
at $x \in \mathbb{V}$, of coefficient $\varepsilon>0$, is the function 
$$\delta^{x}_{\varepsilon}: \mathbb{V} \rightarrow \mathbb{V} \quad , \quad 
\delta^{x}_{\varepsilon} y = x + \varepsilon (-x+y) \quad . $$
For fixed $x$ the dilatations based at $x$ form a one parameter group which
 contracts any bounded neighbourhood of $x$ to a point, uniformly with respect 
to $x$. Thus (Af1'), (A2) are satisfied. (A3') is also obvious.

 The meaning of (A4) is that using dilatations we can recover the operation of addition and multiplication by scalars. 
We shall explain this in detail since this will help the understanding of the axioms of dilatation structures, described in section \ref{secdil}. 
 
For $\displaystyle x,u,v \in \mathbb{V}$ and $\varepsilon>0$ we define the
following compositions of dilatations: 
\begin{equation}
\Delta_{\varepsilon}^{x}(u,v) = \delta_{\varepsilon^{-1}}^{\delta_{\varepsilon}^{x} u}
 \delta^{x}_{\varepsilon} v \quad , 
\label{operations}
\end{equation}
$$\quad 
\Sigma_{\varepsilon}^{x}(u,v) = \delta_{\varepsilon^{-1}}^{x} \delta_{\varepsilon}^{\delta_{\varepsilon}^{x} u}
 (v) \quad , \quad inv^{x}_{\varepsilon}(u) =  \delta_{\varepsilon^{-1}}^{\delta_{\varepsilon}^{x} u} x \quad . $$
The meaning of this functions becomes clear if we compute: 
$$\Delta_{\varepsilon}^{x}(u,v) =  x+\varepsilon(-x+u) + (-u+v)  \quad , $$
$$\Sigma_{\varepsilon}^{x}(u,v) =  u+ \varepsilon(-u+x) + (-x+v) \quad ,$$
$$inv^{x}_{\varepsilon}(u) = =x+\varepsilon(-x+u) + (-u+x) \quad .$$
As $\varepsilon \rightarrow 0$ we have the limits: 
$$\lim_{\varepsilon\rightarrow 0} \Delta_{\varepsilon}^{x}(u,v) = \Delta^{x}(u,v) = x+(-u+v) \quad ,$$
$$\lim_{\varepsilon\rightarrow 0} \Sigma_{\varepsilon}^{x}(u,v) = \Sigma^{x}(u,v) = u+(-x+v) \quad ,$$
$$\lim_{\varepsilon\rightarrow 0} inv^{x}_{\varepsilon}(u) = inv^{x}(u) = x-u+x \quad , $$
uniform with respect to $x,u,v$ in bounded sets. 
The function $\displaystyle  \Sigma^{x}(\cdot,\cdot)$ is a group operation, namely the addition operation 
translated such that the neutral element is $x$: 
$$\Sigma^{x}(u,v) \ = \ u +_{x} v \quad .$$
 The function 
$\displaystyle inv^{x}(\cdot)$ is the inverse function with respect to the operation $\displaystyle +_{x}$
$$inv^{x}(u) \, +_{x} \,u \ = \ u \, +_{x} \, inv^{x}(u) \ = \ x $$ 
 and $\displaystyle  \Delta^{x}(\cdot,\cdot)$ is the difference function
$$\Delta^{x}(u,v) \ = \ inv^{x}(u) \, +_{x}  \, v $$

Notice that for fixed $x, \varepsilon$ the function $\displaystyle  \Sigma^{x}_{\varepsilon}(\cdot,\cdot)$ is not a 
group operation, first of all because it is not associative. Nevertheless, this function satisfies 
a ``shifted'' associativity property, namely 
$$\Sigma_{\varepsilon}^{x}(\Sigma_{\varepsilon}^{x}(u,v),w) = 
\Sigma_{\varepsilon}^{x}(u, \Sigma_{\varepsilon}^{\delta^{x}_{\varepsilon}u}(v,w)) \quad .$$
Also, the inverse function $\displaystyle inv^{x}_{\varepsilon}$ is not involutive, but shifted involutive:  
$$inv_{\varepsilon}^{\delta^{x}_{\varepsilon}u}\left( inv^{x}_{\varepsilon} u\right) = u \quad . $$

Affine continuous  transformations $A:\mathbb{V} \rightarrow \mathbb{V}$ admit the following 
description in terms of dilatations. (We could dispense of continuity hypothesis in this situation, but 
we want to illustrate a general point of view, described further in the paper).  

\begin{prop}
A continuous transformation  $A:\mathbb{V} \rightarrow \mathbb{V}$ is affine if and only if for any 
$\varepsilon \in  (0,1)$, $x,y \in \mathbb{V}$ we have 
\begin{equation}
 A \delta_{\varepsilon}^{x} y \ = \ \delta_{\varepsilon}^{Ax} Ay \quad . 
 \label{eq1proplin} 
 \end{equation}
\label{1proplin}
\end{prop}

The proof is a straightforward consequence of representation formul{\ae} 
(\ref{operations}) for the addition, difference and inverse operations  in terms of dilatations. 

In particular any dilatation is an affine transformation, hence for any 
 $x, y \in \mathbb{V}$ and $\varepsilon, \mu > 0$ we have 
 \begin{equation}
 \delta^{y}_{\mu} \, \delta^{x}_{\varepsilon} \ = \
 \delta^{\delta^{y}_{\mu}x}_{\varepsilon} \delta^{y}_{\mu} \quad . 
 \label{lindil}
 \end{equation}
Thus we recover (Af2') (see also condition (Af2)). The barycentric condition (Af3) is a consequence of the commutativity of the addition of vectors. The easy 
part of the theorem \ref{taffine} is therefore proven. 

The second, difficult part of the theorem is to prove that axioms (Af1'), (A2), (A3'), (A4), 
(Af2') describe normed affine group spaces. This is a direct consequence of several general results from this paper:  theorem \ref{tgrd} and 
proposition \ref{pexlin} show that normed affine group spaces satisfy the axioms, corollary \ref{cortang}, theorem \ref{tdilatlin}, proposition \ref{probari} and theorem \ref{tafflin} show that conversely a space where the axioms are satisfied is a normed affine group space, moreover that in the presence of the 
barycentric condition (Af3) we get real normed affine spaces. \quad $\square$

Some compositions of dilatations are dilatations. This is precisely
stated in the next theorem, which is equivalent with the Menelaos theorem in 
euclidean geometry. 

\begin{thm}
For any  $x, y \in \mathbb{V}$ and $\varepsilon, \mu > 0$ such that $\varepsilon
\mu \not = 1$ there exists an unique
$w \in \mathbb{V}$ such that 
$$\delta^{y}_{\mu} \, \delta^{x}_{\varepsilon} \ = \ 
\delta^{w}_{\varepsilon \mu} \quad . $$
\label{teunu}
\end{thm}

For the proof see Artin \cite{artin}. A straightforward consequence of this theorem is the
following result. 

\begin{cor}
The inverse semigroup  generated by  dilatations of the 
space $\mathbb{V}$ is made of all dilatations and all translations in
$\mathbb{V}$. 
\label{corunu}
\end{cor}

\paragraph{Proof.}
Indeed, by theorem \ref{teunu} a composition of two dilatations with
coefficients $\varepsilon, \mu$ with $\varepsilon \mu \not = 1$ is a dilatation.
By direct computation, if $\varepsilon \mu = 1$ then we obtain
translations. This is in fact compatible with (\ref{operations}), but is a stronger
statement, due to the fact that dilatations are affine in the sense of relation 
(\ref{lindil}).  

Any composition between a translation and a dilatation is again a dilatation. 
The proof is done. \quad $\square$

The corollary \ref{corunu} allows us to describe the ratio of three collinear points in a way which 
will be generalized to normed affine group spaces. Indeed, in a real normed affine space $\mathbb{V}$, for any 
$x, y \in \mathbb{V} $  and $\alpha, \beta \in (0, + \infty)$ such that $\alpha \beta \not = 1$, there is an unique 
$z \in \mathbb{V}$ and $\gamma = 1/ \alpha \beta$ such that 
 $$\delta_{\alpha}^{x} \, \delta^{y}_{\beta} 
\, \delta^{z}_{\gamma} = \, id$$ 
We easily find that $x, y, z$ are collinear 
\begin{equation}
z = \frac{1 - \alpha}{1 - \alpha \beta} x \, + \, \frac{\alpha(1-\beta)}{1
- \alpha \beta}y 
\label{relation}
\end{equation}
the ratio of these three points, named $\displaystyle 
r(x^{\alpha},y^{\beta},z^{\gamma})$ is: 
$$r(x^{\alpha},y^{\beta},z^{\gamma}) \, = \, \frac{\alpha}{1 - \alpha \beta} $$
Conversely, let $x, y, z \in \mathbb{V}$ which are collinear, such that $z$ is in between $x$ and $y$. Then we can easily find (non unique) 
$\alpha, \beta, \gamma \in (0, + \infty)$ such that $\alpha \beta \gamma = 1$ and  $\displaystyle \delta_{\alpha}^{x} \, \delta^{y}_{\beta} 
\, \delta^{z}_{\gamma} = \, id$. 

It is then  almost straightforward to prove the well known fact that any affine transformation is also geometrically affine, in the sense that 
it transforms triples of collinear points into triples of collinear points (use commutation with dilatations) and it preserves the ratio of collinear 
points. (The converse is also true). 

\section{Dilatation structures}
\label{secdil}
A dilatation structure $(X,d,\delta)$ over a metric space $(X,d)$ is an assignment to any 
point $x \in X$ of a group of "dilatations" $\displaystyle \left\{ \delta^{x}_{\varepsilon} \mbox{ : } 
\varepsilon \in \Gamma \right\}$, together with some compatibility conditions between the distance 
and the dilatations and between dilatations based in different points. 

A basic difficulty in stating the axioms of a dilatation structure is related to the domain of definition and 
the image of a dilatation. In this subsection we shall neglect the problems raised by domains and codomains of 
dilatations. 

The axioms state that some combinations between dilatations and the distance converge uniformly, with respect 
to some finite families of points in an arbitrary compact subset of the metric space $(X,d)$, as $\nu(\varepsilon)$ 
converges to $0$.

We present here an  introduction into the subject of dilatation 
structures. For more details see Buliga \cite{buligadil1}.

\subsection{Notations}

Let $\Gamma$ be  a topological separated commutative group  endowed with a continuous group morphism 
$\nu : \Gamma \rightarrow (0,+\infty)$ with $\displaystyle \inf \nu(\Gamma)  =  0$. Here $(0,+\infty)$ is 
taken as a group with multiplication. The neutral element of $\Gamma$ is denoted by $1$. We use the multiplicative notation for the operation in $\Gamma$. 

The morphism $\nu$ defines an invariant topological filter on $\Gamma$ (equivalently, an end). Indeed, 
this is the filter generated by the open sets $\displaystyle \nu^{-1}(0,a)$, $a>0$. From now on 
we shall name this topological filter (end) by "0" and we shall write $\varepsilon \in \Gamma \rightarrow 
0$ for $\nu(\varepsilon)\in (0,+\infty) \rightarrow 0$. 

The set $\displaystyle \Gamma_{1} = \nu^{-1}(0,1] $ is a semigroup. We note $\displaystyle 
\bar{\Gamma}_{1}= \Gamma_{1} \cup \left\{ 0\right\}$
On the set $\displaystyle 
\bar{\Gamma}= \Gamma \cup \left\{ 0\right\}$ we extend the operation on $\Gamma$ by adding the rules  
$00=0$ and $\varepsilon 0 = 0$ for any $\varepsilon \in \Gamma$. This is in agreement with the invariance 
of the end $0$ with respect to translations in $\Gamma$.

The space $(X,d)$ is a complete, locally compact metric space. For any $r>0$  
and any $x \in X$ we denote by $B(x,r)$ the open ball of center $x$ and radius 
$r$ in the metric space $X$.

On the metric space $(X,d)$ we work with the topology (and uniformity) induced 
by the distance. For any $x \in X$ we denote by $\mathcal{V}(x)$ the topological
filter of open neighbourhoods of $x$. 

The dilatation structures, which will be introduced soon, are invariant to the operation of multiplication 
of the distance by a positive constant. They should also be seen, as examples show, as local objects, therefore 
we may safely suppose, without restricting the generality, that all closed balls of radius at most $5$ are compact.

\subsection{Axioms of dilatation structures}

We shall list the   axioms of  a dilatation structure $(X,d,\delta)$, in a simplified form, without concerning about 
domains and codomains of functions. In the next subsection we shall add the supplementary conditions concerning domains and codomains 
of dilatations.

\begin{enumerate}
\item[{\bf A1.}] {\it For any point $x \in X$ there is an action $\displaystyle \delta^{x}: \Gamma \rightarrow 
End(X,d, x)$, where $End(X,d, x)$ is the collection of all continuous, with continuous inverse transformations  
$\phi: (X,d) \rightarrow (X,d)$ such that $\phi(x) = x$.} 
\end{enumerate}

This axiom (the same as (A1) from theorem \ref{tbertram} or theorem \ref{taffine}) tells us that $\displaystyle \delta^{x}_{\varepsilon} x = x$ for any $x \in X$, $\varepsilon \in \Gamma$, 
also $\displaystyle \delta^{x}_{1} y = y$ for any $x,y \in X$, and $\displaystyle \delta^{x}_{\varepsilon} \delta^{x}_{\mu}
 y = \delta^{x}_{\varepsilon \mu} y$ for any $x, y \in X$ and $\varepsilon, \mu \in \Gamma$.

\begin{enumerate}
\item[{\bf A2.}]  
{\it The function $\displaystyle \delta : \Gamma \times X \times X \rightarrow  X$ defined by 
$\displaystyle \delta (\varepsilon,  x, y)  = \delta^{x}_{\varepsilon} y$ is continuous. Moreover, it can be 
continuously extended to $\displaystyle \bar{\Gamma} \times X \times X$ by $\delta (0,x, y) = x$  and the limit  
$$\lim_{\varepsilon\rightarrow 0} \delta_{\varepsilon}^{x} y \, = \, x  $$
is uniform with respect to $x,y$ in compact set.}
\end{enumerate}

We may alternatively put that the previous limit is uniform with respect to $d(x,y)$. 

\begin{enumerate}
\item[{\bf A3.}] {\it There is $A > 1$ such that  for any $x$ there exists 
 a  function $\displaystyle (u,v) \mapsto d^{x}(u,v)$, defined for any 
$u,v$ in the closed ball (in distance d) $\displaystyle 
\bar{B}(x,A)$, such that 
$$\lim_{\varepsilon \rightarrow 0} \quad \sup  \left\{  \mid \frac{1}{\varepsilon} d(\delta^{x}_{\varepsilon} u, 
\delta^{x}_{\varepsilon} v) \ - \ d^{x}(u,v) \mid \mbox{ :  } u,v \in \bar{B}_{d}(x,A)\right\} \ =  \ 0$$
uniformly with respect to $x$ in compact set. }
\end{enumerate}

It is easy to see that: 
\begin{enumerate}
\item[(a)] The function  $\displaystyle d^{x}$ is continuous as an uniform limit of continuous functions on a compact set, 
\item[(b)] $\displaystyle d^{x}$ is symmetric $\displaystyle d^{x}(u,v) = d^{x}(v,u)$ for any $u,v \in \bar{B}(x,A)$, 
\item[(c)] $\displaystyle d^{x}$ satisfies the triangle inequality, but it can be a degenerated distance function: 
there might exist  $\displaystyle v,w $ such that $\displaystyle d^{x}(v,w) = 0$. 
\end{enumerate}

We make the following notation which generalizes  the notation from (\ref{operations}): 
$$\Delta^{x}_{\varepsilon}(u,v) = \delta_{\varepsilon^{-1}}^{\delta^{x}_{\varepsilon} u} \delta^{x}_{\varepsilon} v . $$
The next axiom can now be stated: 
\begin{enumerate}
\item[{\bf A4.}] We have the limit 
$$\lim_{\varepsilon \rightarrow 0}  \Delta^{x}_{\varepsilon}(u,v) =  \Delta^{x}(u, v)  $$
uniformly with respect to $x, u, v$ in compact set. 
\end{enumerate}

\begin{dfn}
A triple $(X,d,\delta)$ which satisfies  A1, A2, A3, but $\displaystyle d^{x}$ is degenerate for some 
$x\in X$, is called degenerate dilatation structure. 

If the triple $(X,d,\delta)$ satisfies  A1, A2, A3 and 
 $\displaystyle d^{x}$ is non-degenerate for any $x\in X$, then we call it  a 
 dilatation structure. 

 If a  dilatation structure satisfies A4 then we call it strong dilatation 
 structure. 
 \label{defweakstrong}
\end{dfn}

\subsection{Axiom $0$: domains and codomains of dilatations}

The problem of domains and codomains of dilatation cannot be neglected. In the section dedicated to examples of dilatation structures  
we present the particular case of an ultrametric space which is also a ball of radius one. As dilatations approximately contract distances, 
it follows that the codomain of a dilatation $\displaystyle \delta^{x}_{\varepsilon}$ with $\nu(\varepsilon) < 1$ can not be the whole space. 
There are other examples showing that we can not always take the domain of a dilatation to be the whole space. That is because the topology of 
small balls can be different from the topology of big ones (like in the case of  a sphere). 

For all these reasons we need to impose some minimal conditions on the domains and codomains of dilatations. These conditions will be explained 
in the following. They will be considered as part of a new axiom, called Axiom $0$. 

For any $x \in X$ there is an open neighbourhood $U(x)$ of $x$ such that for any $\displaystyle \varepsilon \in \Gamma_{1}$ 
the dilatations are functions  $$ \delta_{\varepsilon}^{x}: U(x) \rightarrow V_{\varepsilon}(x) \quad . $$ 
 The sets $\displaystyle  V_{\varepsilon}(x)$ are open neighbourhoods of $x$.  

There is a number  $1<A$ such that for any $x \in X$ we have $\displaystyle \bar{B}_{d}(x,A) \subset U(x)$. 
There is a number $B > A$ such that  for 
 any $\varepsilon \in \Gamma$ with $\nu(\varepsilon) \in (1,+\infty)$ the 
 associated dilatation  is a function 
$$\delta^{x}_{\varepsilon} : W_{\varepsilon}(x) \rightarrow B_{d}(x,B) \quad .$$

We have   the following string of inclusions, for any $\displaystyle \varepsilon \in \Gamma_{1}$, and any $x \in X$:
$$ B_{d}(x,\nu(\varepsilon)) \subset \delta^{x}_{\varepsilon}  B_{d}(x, A) 
\subset V_{\varepsilon}(x) \subset 
W_{\varepsilon^{-1}}(x) \subset \delta_{\varepsilon}^{x}  B_{d}(x, B) \quad . $$

In relation with the axiom A4 we need the  following  condition on the co-domains $\displaystyle V_{\varepsilon}(x)$: for any compact set $K \subset X$ there are $R=R(K) > 0$ and $\displaystyle \varepsilon_{0}= \varepsilon(K) \in (0,1)$  such that  
for all $\displaystyle u,v \in \bar{B}_{d}(x,R)$ and all $\displaystyle \varepsilon \in \Gamma$, $\displaystyle  \nu(\varepsilon) \in (0,\varepsilon_{0})$,  we have 
$$\delta_{\varepsilon}^{x} v \in W_{\varepsilon^{-1}}( \delta^{x}_{\varepsilon}u) \ .$$

These conditions are important for describing dilatation structures on the boundary of the dyadic tree, for example. In the first formulation of the axioms given in  \cite{buligadil1} some of these assumptions are part of the Axiom $0$, others can be found in the initial formulation of the Axioms 1, 2, 3.

\section{Groups with dilatations}

For a dilatation structure the metric tangent spaces   have a group structure which is compatible 
with dilatations. This structure, of a normed group with dilatations, is interesting 
by itself. The notion has been introduced in \cite{buliga2}, \cite{buligadil1}; we describe  it 
further.

 We shall work further with local groups. Such objects are not groups:  they are spaces endowed with an operation 
defined only locally, satisfying the conditions of a uniform group. In \cite{buliga2} we use a slightly non standard 
definition of such objects. For the purposes of this paper it seems enough to mention that neighbourhoods of the 
neutral element in a uniform group are local groups. See section 3.3 \cite{buligadil1} for details about the definition 
of local groups.

\begin{dfn}
A group with dilatations $(G,\delta)$ is a local  group $G$  with  a local action of $\Gamma$ (denoted by $\delta$), on $G$ such that
\begin{enumerate}
\item[H0.] the limit  $\displaystyle \lim_{\varepsilon \rightarrow 0} 
\delta_{\varepsilon} x  =  e$ exists and is uniform with respect to $x$ in a compact neighbourhood of the identity $e$.
\item[H1.] the limit
$$\beta(x,y)  =  \lim_{\varepsilon \rightarrow 0} \delta_{\varepsilon}^{-1}
\left((\delta_{\varepsilon}x) (\delta_{\varepsilon}y ) \right)$$
is well defined in a compact neighbourhood of $e$ and the limit is uniform.
\item[H2.] the following relation holds
$$ \lim_{\varepsilon \rightarrow 0} \delta_{\varepsilon}^{-1}
\left( ( \delta_{\varepsilon}x)^{-1}\right)  =  x^{-1}$$
where the limit from the left hand side exists in a neighbourhood of $e$ and is uniform with respect to $x$.
\end{enumerate}
\label{defgwd}
\end{dfn}

\begin{dfn} A normed group with dilatations $(G, \delta, \| \cdot \|)$ is a 
group with dilatations  $(G, \delta)$ endowed with a continuous norm  
function $\displaystyle \|\cdot \| : G \rightarrow \mathbb{R}$ which satisfies 
(locally, in a neighbourhood  of the neutral element $e$) the properties: 
 \begin{enumerate}
 \item[(a)] for any $x$ we have $\| x\| \geq 0$; if $\| x\| = 0$ then $x=e$, 
 \item[(b)] for any $x,y$ we have $\|xy\| \leq \|x\| + \|y\|$, 
 \item[(c)] for any $x$ we have $\displaystyle \| x^{-1}\| = \|x\|$, 
 \item[(d)] the limit 
$\displaystyle \lim_{\varepsilon \rightarrow 0} \frac{1}{\nu(\varepsilon)} \| \delta_{\varepsilon} x \| = \| x\|^{N}$ 
 exists, is uniform with respect to $x$ in compact set, 
 \item[(e)] if $\displaystyle \| x\|^{N} = 0$ then $x=e$.
  \end{enumerate}
  \label{dnco}
  \end{dfn}

In a normed group with dilatations we have a natural left invariant distance given by
\begin{equation}
d(x,y) = \| x^{-1}y\| \quad . 
\label{dnormed}
\end{equation}
Any normed group with dilatations has an associated dilatation structure on it.  In a group with dilatations $(G, \delta)$  we define dilatations based in any point $x \in G$ by 
 \begin{equation}
 \delta^{x}_{\varepsilon} u = x \delta_{\varepsilon} ( x^{-1}u)  . 
 \label{dilat}
 \end{equation}

The following result is theorem 15 \cite{buligadil1}. 

\begin{thm}
Let $(G, \delta, \| \cdot \|)$ be  a locally compact  normed local group with dilatations. Then $(G, d, \delta)$ is 
a dilatation structure, where $\delta$ are the dilatations defined by (\ref{dilat}) and the distance $d$ is induced by the norm as in (\ref{dnormed}). 
\label{tgrd}
\end{thm}

\subsection{Conical groups}

\begin{dfn}
A normed conical group $N$ is a normed  group with dilatations  such that for any $\varepsilon \in \Gamma$  the dilatation 
 $\delta_{\varepsilon}$ is a group morphism  and such that for any $\varepsilon >0$ 
  $\displaystyle  \| \delta_{\varepsilon} x \| = \nu(\varepsilon) \| x \|$. 
\end{dfn}

A conical group is the infinitesimal version of a group with 
dilatations (\cite{buligadil1} proposition 2).

\begin{prop}
Under the hypotheses H0, H1, H2 $\displaystyle (G,\beta, \delta, \| \cdot \|^{N})$ is a local normed conical group, with operation 
$\beta$,  dilatations $\delta$ and homogeneous norm $\displaystyle \| \cdot \|^{N}$.
\label{here3.4}
\end{prop}

\subsection{Carnot groups}
\label{carnotgroups}

Carnot groups appear in sub-riemannian geometry 
 as models of tangent spaces,   \cite{bell}, \cite{gromovgr}, \cite{pansu}. In particular such groups can be endowed with a structure 
of sub-riemannian manifold.

\begin{dfn}
A Carnot (or stratified homogeneous) group is a pair $\displaystyle 
(N, V_{1})$ consisting of a real 
connected simply connected group $N$  with  a distinguished subspace  
$V_{1}$ of  the Lie algebra $Lie(N)$, such that  the following   
direct sum decomposition occurs: 
$$n \ = \ \sum_{i=1}^{m} V_{i} \ , \ \ V_{i+1} \ = \ [V_{1},V_{i}]$$
The number $m$ is the step of the group. The number $\displaystyle Q \ = \ \sum_{i=1}^{m} i 
\ dim V_{i}$ is called the homogeneous dimension of the group. 
\label{dccgroup}
\end{dfn}

Because the group is nilpotent and simply connected, the
exponential mapping is a diffeomorphism. We shall identify the 
group with the algebra, if is not locally otherwise stated.

The structure that we obtain is a set $N$ endowed with a Lie
bracket and a group multiplication operation, related by the 
Baker-Campbell-Hausdorff formula. Remark that  the group operation
is polynomial.

Any Carnot group admits a one-parameter family of dilatations. For any 
$\varepsilon > 0$, the associated dilatation is: 
$$ x \ = \ \sum_{i=1}^{m} x_{i} \ \mapsto \ \delta_{\varepsilon} x \
= \ \sum_{i=1}^{m} \varepsilon^{i} x_{i}$$
Any such dilatation is a group morphism and a Lie algebra morphism.

In a Carnot group $N$ let us choose  an  euclidean norm 
$\| \cdot \|$ on $\displaystyle V_{1}$.  We shall endow the group $N$ with 
a structure of a sub-riemannian manifold. For this take the distribution 
obtained from left translates of the space $V_{1}$. The metric on that 
distribution is obtained by left translation of the inner product restricted to
$V_{1}$.

Because $V_{1}$  generates (the algebra) $N$
then any element $x \in N$ can be written as a product of 
elements from $V_{1}$, in a controlled way, described in the following  useful lemma  (slight reformulation of Lemma 1.40, 
Folland, Stein \cite{fostein}).  

\begin{lem}
Let $N$ be a Carnot group and $X_{1}, ..., X_{p}$ an orthonormal basis 
for $V_{1}$. Then there is a  
 a natural number $M$ and a function $g: \left\{ 1,...,M \right\} 
\rightarrow \left\{ 1,...,p\right\}$ such that any element 
$x \in N$ can be written as: 
\begin{equation}
x \ = \ \prod_{i = 1}^{M} \exp(t_{i}X_{g(i)})
\label{fp2.4}
\end{equation}
Moreover, if $x$ is sufficiently close (in Euclidean norm) to
$0$ then each $t_{i}$ can be chosen such that $\mid t_{i}\mid \leq C 
\| x \|^{1/m}$
\label{p2.4}
\end{lem}

As a consequence we get: 

\begin{cor}
 The Carnot-Carath\'eodory distance 
$$d(x,y) \ = \ \inf \left\{ \int_{0}^{1} \| c^{-1} \dot{c} \| \mbox{ 
d}t \ \mbox{ : } \ c(0) = x , \ c(1) = y , \quad \quad \right.$$ 
$$\left. \quad \quad \quad \quad  c^{-1}(t) \dot{c}(t) \in V_{1} 
\mbox{ for a.e. } t \in [0,1] 
\right\}$$ 
is finite for any two $x,y \in N$.  The distance is obviously left
invariant, thus it induces a norm on $N$. 
\end{cor}

The Carnot-Carath\'eodory distance induces a homogeneous norm on the Carnot group $N$ by the formula: 
$\| x \| = d(0,x)$. From the invariance of the distance with respect to left translations we get: for any $x, y \in N$  
$$\|x^{-1} y\| = d(x,y)$$

For any $x \in N$ and $\varepsilon > 0$ we define the dilatation $\displaystyle \delta^{x}_{\varepsilon} y = x \delta_{\varepsilon}(x^{-1} y)$. 
Then $(N, d, \delta)$ is a dilatation structure, according to theorem \ref{tgrd}.

\subsection{Contractible groups}

\begin{dfn}
A contractible group is a pair $(G,\alpha)$, where $G$ is a  
topological group with neutral element denoted by $e$, and $\alpha \in Aut(G)$ 
is an automorphism of $G$ such that: 
\begin{enumerate}
\item[-] $\alpha$ is continuous, with continuous inverse, 
\item[-] for any $x \in G$ we have the limit $\displaystyle 
\lim_{n \rightarrow \infty} \alpha^{n}(x) = e$. 
\end{enumerate}
\end{dfn}

For a contractible group  $(G,\alpha)$, the automorphism $\alpha$ is compactly contractive (Lemma 1.4 (iv) \cite{siebert}), that is: 
for each compact set 
$K \subset G$ and open set $U \subset G$, with $e \in U$, there is 
$\displaystyle N(K,U) \in \mathbb{N}$ such that for any $x \in K$ and $n \in
\mathbb{N}$, $n \geq N(K,U)$, we have $\displaystyle \alpha^{n}(x) \in U$. 

If $G$ is locally compact then $\alpha$ compactly contractive is equivalent with:  each identity neighbourhood 
of $G$ contains an $\alpha$-invariant neighbourhood. Further on we shall assume without mentioning 
that all groups are locally compact.

Any conical group can be seen as a contractible  group. 
Indeed, it suffices to associate to a conical group $(G,\delta)$ the 
contractible group $\displaystyle (G,\delta_{\varepsilon})$, for a fixed 
$\varepsilon \in \Gamma$ with $\nu(\varepsilon) < 1$.

Conversely, to any contractible group $(G,\alpha)$ we may  associate the conical group $(G,\delta)$, with 
$\displaystyle \Gamma = \left\{ \frac{1}{2^{n}} \mbox{ : } n \in \mathbb{N}
\right\}$ and for any $n \in \mathbb{N}$ and $x \in G$ 
$$\displaystyle \delta_{\frac{1}{2^{n}}} x \ = \ \alpha^{n}(x) \quad . $$
(Finally, a local conical
group has only locally the structure of a contractible group.)

The structure of contractible groups is known in some detail, due to Siebert
\cite{siebert}, Wang \cite{wang}, Gl\"{o}ckner and Willis \cite{glockwill}, 
Gl\"{o}ckner \cite{glockner} and references therein. 

For this paper the following results are of interest. We begin with  the
definition of a contracting automorphism group \cite{siebert}, definition 5.1. 

\begin{dfn}
Let $G$ be a locally compact group. An automorphism group on $G$ is a family 
$\displaystyle T= \left( \tau_{t}\right)_{t >0}$ in $Aut(G)$, such that 
$\displaystyle \tau_{t} \, \tau_{s} = \tau_{ts}$ for all $t,s > 0$. 

The contraction group of $T$ is defined by 
$$C(T) \ = \ \left\{ x \in G \mbox{ : } \lim_{t \rightarrow 0} \tau_{t}(x) = e
\right\} \quad .$$
The automorphism group $T$ is contractive if $C(T) = G$. 
\end{dfn}

It is obvious that a contractive automorphism  group $T$ induces on $G$ a 
structure of conical group. Conversely, any conical group with $\Gamma = 
(0,+\infty)$ has an associated contractive automorphism group (the group of 
dilatations based at the neutral element). 

Further is proposition 5.4 \cite{siebert}. 

\begin{prop}
For a locally compact group $G$ the following assertions are equivalent: 
\begin{enumerate}
\item[(i)] $G$ admits a contractive automorphism group;
\item[(ii)] $G$ is a simply connected Lie group whose Lie algebra admits a 
positive graduation.
\end{enumerate}
\label{psiebert}
\end{prop}

\section{Other examples of dilatation structures}

\subsection{Riemannian manifolds} 

The following interesting   quotation from Gromov book \cite{gromovbook}, pages 85-86, motivates some of the  ideas underlying  dilatation structures, 
especially in the very particular case of a riemannian manifold: 

``{\bf 3.15. 
Proposition:} {\it Let $(V, g)$ be a Riemannian manifold with $g$ continuous.
For each $v \in V$ the spaces $(V, \lambda d , v)$ Lipschitz converge as  $\lambda \rightarrow \infty$ to the 
tangent space $(T_{v}V, 0)$ with its Euclidean metric $g_{v}$.} 

$\mathbf{Proof_{+}:}$ {\it Start with a $C^{1}$ map $\mathbb{(R}^{n}, 0) \rightarrow  (V, v)$ whose differential is 
isometric at 0. The $\lambda$-scalings of this provide almost isometries between large balls in $\mathbb{R}^{n}$
 and those in $\lambda V$ for  $\lambda \rightarrow \infty$.} 
 {\bf Remark:} {\it In fact we can define Riemannian manifolds as 
 locally compact path metric spaces that satisfy the conclusion of Proposition 3.15.}``

The problem of domains and codomains left aside, any chart of a  Riemannian manifold induces locally a dilatation 
structure on the manifold. Indeed, take $(M,d)$ to be a $n$-dimensional  Riemannian manifold with $d$ the distance on 
$M$ induced by the Riemannian structure. Consider a diffeomorphism $\phi$ of an open set $U \subset M$ onto 
$\displaystyle V \subset \mathbb{R}^{n}$  and transport the dilatations from $V$ to $U$ (equivalently, transport the distance $d$ from 
$U$ to $V$). There is only one thing to check in order to see that we got a dilatation structure: the axiom A3, expressing 
the compatibility of the distance $d$ with the dilatations. But this is just a metric way to express the distance 
on the tangent space of $M$ at $x$ as a limit of rescaled distances (see  Gromov Proposition 3.15, \cite{gromovbook},
p. 85-86). Denoting by $\displaystyle g_{x}$ the metric tensor at $x \in U$, we have:  
$$\left[ d^{x}(u,v) \right]^{2}  = $$ 
$$= g_{x}\left( \frac{d}{d \, \varepsilon}_{|_{\varepsilon = 0}}\phi^{-1}\left(\phi(x) + \varepsilon (\phi(u) -
\phi(x))\right) ,  \frac{d}{d \, \varepsilon}_{|_{\varepsilon = 0}}\phi^{-1}\left(\phi(x) + \varepsilon (\phi(v) -
\phi(x))\right) \right)$$

A basically different  example of a dilatation structure on a riemannian manifold  will be explained next. Let $M$ be 
a $n$ dimensional riemannian manifold and $\exp$ be the geodesic exponential. To any point $x \in M$ and any vector 
$\displaystyle v \in T_{x} M$ the point   $\displaystyle \exp_{x}(v) \in M$ is located on the geodesic passing thru $x$ 
and tangent to $v$; if we parameterize this geodesic with respect to length, such that the tangent at $x$ is parallel and 
has the same direction as $v$, then $\displaystyle \exp_{x}(v) \in M$ has the coordinate equal with the length of $v$ with 
respect to the norm on $\displaystyle T_{x} M$.  We define implicitly  the dilatation based at $x$, of coefficient 
$\varepsilon > 0$ by the relation: 
 $$\delta^{x}_{\varepsilon} \exp_{x} (u) \, = \, \exp_{x} \left( \varepsilon u \right) \quad . $$
It is not straightforward to check that we obtain a strong dilatation structure, but it is true. There are interesting facts 
related to the numbers $A, B$ and  the minimal regularity required for the riemannian manifold. This example is different from the 
first because instead of using a chart (same for all $x$) we use a family of charts indexed with respect to the basepoint of the 
dilatations.

\subsection{Dilatation structures on the boundary of the dyadic tree}

We shall take the group $\Gamma$ to be the set of integer powers of $2$, seen as a subset of 
dyadic numbers. Thus for any $p \in \mathbb{Z}$ the element $\displaystyle 2^{p} \in \mathbb{Q}_{2}$ 
belongs to $\Gamma$. The operation is the multiplication of dyadic numbers and the morphism 
$\nu : \Gamma \rightarrow (0,+\infty)$ is defined by 
$$\nu(2^{p}) = d(0, 2^{p}) = \frac{1}{2^{p}} \in (0,+\infty) \quad  . $$

The dyadic tree $\mathcal{T}$  is the infinite rooted planar binary tree.  
Any node has two descendants. The nodes are  coded by elements of 
$\displaystyle X^{*}$, $X = \left\{ 0,1\right\}$. The root is coded by the 
empty word and if a node is coded by $x\in X^{*}$ then its left hand side 
descendant has the code $x0$ and its 
right hand side descendant has the code $x1$. We shall therefore identify the 
dyadic tree with $\displaystyle X^{*}$ and we put on the dyadic tree the 
natural (ultrametric) distance on $\displaystyle X^{*}$. The boundary (or 
the set of ends) of the dyadic tree is then the same as the compact 
ultrametric space $\displaystyle X^{\omega}$. 

 $\displaystyle X^{\omega}$ is the set of 
words infinite at right over the alphabet $X = \left\{ 0, 1 \right\}$: 
$$X^{\omega} = \left\{ f \  \mid \ \   f: \mathbb{N}^{*}\rightarrow X \right\} = X^{\mathbb{N}^{*}} \quad . $$
A natural distance on this set is defined for different $\displaystyle x , y \in X^{\omega}$ by the formula
$$d(x,y) = \frac{1}{2^{m}}$$ 
where $m$ is the length of largest common prefix of the words $x$ and $y$. This distance is ultrametric. The metric space 
$\displaystyle (X^{\omega}, d)$ is isometric with the space of dyadic integers. The metric space is then a ball of radius $1$. 

A trivial dilatation structure is induced by the identification with dyadic integers and it has the following expression: 
$$\delta^{x}_{2^{p}} y \, = \, x + 2^{p} ( y - x) $$
where the operations are done with dyadic integers. 

More complex dilatation structures are given by the following construction. See theorem 6.5 \cite{buligaself} for more details.

\begin{dfn}
A function $\displaystyle W: \mathbb{N}^{*} \times X^{\omega} \rightarrow 
Isom(X^{\omega})$  is smooth if for any $\varepsilon > 0$ there exists $\mu(\varepsilon) > 0$ such that 
for any $\displaystyle x, x' \in X^{\omega}$ such that $d(x,x')< \mu(\varepsilon)$ and for any $\displaystyle y \in X^{\omega}$  we have 
$$ \frac{1}{2^{k}} \, d( W^{x}_{k} (y) , W^{x'}_{k} (y) ) \leq \varepsilon \quad , $$
for an  $k$ such that  $\displaystyle d(x,x') <   1 / 2^{k} $.  
\label{defwsmooth}
\end{dfn}

\begin{thm}
To any smooth function  $\displaystyle W: \mathbb{N}^{*} \times X^{\omega} \rightarrow 
Isom(X^{\omega})$ in the sense of definition \ref{defwsmooth}  is associated a  dilatation structure  $\displaystyle (X^{\omega}, d, \delta)$, 
 induced by  functions $\displaystyle \delta_{2}^{x}$, defined by 
$\displaystyle \delta_{2}^{x} x = x$ and otherwise by:  for any 
$\displaystyle q \in X^{*}$, $\alpha \in X$, $x, y \in X^{\omega}$ we have 
\begin{equation}
\delta_{2}^{q \alpha x} q \bar{\alpha} y = q \alpha \bar{x_{1}} W^{q \alpha x}_{\mid q \mid + 1} 
(y) \quad  . 
\label{eqtstruc}
\end{equation}
\end{thm}

\subsection{Sub-riemannian manifolds}

Regular sub-riemannian manifolds provide examples of dilatation structures. In the paper \cite{buligasr} this is explained in all details. 
See section \ref{carnotgroups} for the most basic example of a dilatation structure on a sub-riemannian manifold: the case of a Carnot group. 

More general,  the dilatation structures constructed over normed groups with dilatations (theorem \ref{tgrd}), 
with $\Gamma = (0,+\infty)$ and $\nu = \, id$, provide more examples of sub-riemannian dilatation structures. 

A sub-riemannian  manifold is a triple $(M,D, g)$, where $M$ is a
connected manifold, $D$ is a completely non-integrable distribution on $M$, and $g$ is a metric (Euclidean inner-product) on the distribution 
(or horizontal bundle)  $D$. A horizontal curve $c:[a,b] \rightarrow M$ is a curve which is almost everywhere derivable and for
almost any $t \in [a,b]$ we have $\displaystyle \dot{c}(t) \in D_{c(t)}$. The class of horizontal curves is denoted by 
$Hor(M,D)$. The following theorem of Chow \cite{chow} is well known.

\begin{thm} (Chow) Let $D$ be a distribution of dimension $m$  in the manifold
$M$. Suppose there is a positive integer number $k$ (called the rank of the
distribution $D$) such that for any $x \in X$ there is a topological  open ball
$U(x) \subset M$ with $x \in U(x)$ such that there are smooth vector fields
$\displaystyle X_{1}, ..., X_{m}$ in $U(x)$ with the property:

(C) the vector fields $\displaystyle X_{1}, ..., X_{m}$ span $\displaystyle
D_{x}$ and these vector fields together with  their iterated
brackets of order at most $k$ span the tangent space $\displaystyle T_{y}M$
at every point $y \in U(x)$.

Then $M$ is locally connected
by horizontal curves 
\label{tchow}
\end{thm}

The Carnot-Carath\'eodory distance (or CC distance) associated to the sub-riemannian manifold is the
distance induced by the length $l$ of horizontal curves:
$$d(x,y) \ = \ \inf \left\{ l(c) \mbox{ : } c \in Hor(M,D) \
, \ c(a) = x \ , \  c(b) = y \right\} $$

Chow condition (C) is used to construct an adapted frame starting from a family of vector fields which generate 
the distribution $D$. A fundamental result in sub-riemannian geometry is the existence of normal frames. This existence 
result is based on the accumulation of various results by Bella\"{\i}che \cite{bell}, first to speak about normal frames, providing rigorous proofs 
for this existence in a flow of results between theorem 4.15 and ending in the first half of section 7.3 (page 62), Gromov \cite{gromovsr} in his 
approximation theorem  p. 135 (conclusion of the point (a) below), as well in his convergence results concerning the nilpotentization of vector fields 
(related to point (b) below), Vodopyanov and others \cite{vodopis} \cite{vodopis2} 
 \cite{vodokar} concerning the proof of basic results in sub-riemannian geometry under very weak regularity assumptions 
(for a discussion of this see \cite{buligasr}). There is no place here to submerge into this, we shall just assume that the object defined below exists.

\begin{dfn}
An adapted frame $\displaystyle \left\{ X_{1}, ... , X_{n} \right\}$ is a normal
frame if the following two conditions are satisfied:
\begin{enumerate}
\item[(a)] we have the limit
$$\lim_{\varepsilon \rightarrow 0_{+}} \frac{1}{\varepsilon} \, d\left(
\exp \left( \sum_{1}^{n}\varepsilon^{deg\, X_{i}} a_{i} X_{i} \right)(y), y \right) \ = \ A(y, a) \in
(0,+\infty)$$
uniformly with respect to $y$ in compact sets and $\displaystyle a=(a_{1}, ...,
a_{n}) \in W$, with $\displaystyle W \subset \mathbb{R}^{n}$ compact
neighbourhood of $\displaystyle 0 \in \mathbb{R}^{n}$,
\item[(b)] for any compact set $K\subset M$ with diameter (with respect to the
distance $d$) sufficiently small,  and for any $i = 1, ..., n$ there
are functions 
$$ P_{i}(\cdot, \cdot, \cdot): U_{K} \times U_{K} \times K \rightarrow  \mathbb{R}$$
 with $\displaystyle U_{K} \subset \mathbb{R}^{n}$ a sufficiently small compact neighbourhood of 
$\displaystyle 0 \in \mathbb{R}^{n}$ such that for any   $x \in K$ and any $\displaystyle 
a,b \in U_{K}$ we have
$$\exp \left( \sum_{1}^{n} a_{i} X_{i} \right) (x) \ = \
\exp \left( \sum_{1}^{n} P_{i}(a, b, y) X_{i} \right) \circ \exp \left( \sum_{1}^{n}  b_{i} X_{i} \right) (x) $$
and such that the following limit exists
$$\lim_{\varepsilon \rightarrow 0_{+}}
\varepsilon^{-deg \, X_{i}} P_{i}(\varepsilon^{deg \, X_{j}} a_{j}, \varepsilon^{deg \, X_{k}} b_{k}, x)   \in
\mathbb{R}$$
and it is uniform with respect to $x  \in K$ and $\displaystyle a, b \in U_{K}$.
\end{enumerate}
\label{defnormal}
\end{dfn}

With the help of a normal frame we can prove the existence of strong dilatation structures on regular sub-riemannian manifolds. The following is 
a consequence of theorems 6.3, 6.4 \cite{buligasr}. 

\begin{thm} 
 Let $(M,D, g)$ be a regular sub-riemannian manifold, $U \subset M$ an open set which admits a normal frame. Define for any $x \in  U$  
and $\varepsilon > 0$ (sufficiently small if necessary),
the dilatation $\displaystyle \delta^{x}_{\varepsilon}$  given by:
$$\delta^{x}_{\varepsilon} \left(\exp\left( \sum_{i=1}^{n} a_{i} X_{i} \right)(x)\right) \  = \
\exp\left( \sum_{i=1}^{n} a_{i} \varepsilon^{deg X_{i}}  X_{i} \right)(x)$$
Then $(U, d,  \delta)$ is a strong dilatation structure. 
\end{thm}

\section{Properties of dilatation structures}

\subsection{Tangent bundle}

A reformulation of parts of theorems 6,7 \cite{buligadil1} is the following. 

\begin{thm}
A   dilatation structure $(X,d,\delta)$ has the following properties.  
\begin{enumerate}
\item[(a)] For all $x\in X$, $u,v \in X$ such that $\displaystyle d(x,u)\leq 1$ and $\displaystyle d(x,v) \leq 1$  and all $\mu \in (0,A)$ we have: 
$$d^{x}(u,v) \ = \ \frac{1}{\mu} d^{x}(\delta_{\mu}^{x} u , \delta^{x}_{\mu} v) \ .$$
We shall say that $d^{x}$ has the cone property with respect to dilatations. 
\item[(b)] The metric space $(X,d)$ admits a metric tangent space 
at $x$, for any point $x\in X$. More precisely we have  the following limit: 
$$\lim_{\varepsilon \rightarrow 0} \ \frac{1}{\varepsilon} \sup \left\{  \mid d(u,v) - d^{x}(u,v) \mid \mbox{ : } d(x,u) \leq \varepsilon \ , \ d(x,v) \leq \varepsilon \right\} \ = \ 0 \ .$$
\end{enumerate}
\label{thcone}
\end{thm}

For the next theorem (composite of results in theorems 8, 10 \cite{buligadil1}) 
 we need the previously introduced notion  of a normed conical local group.

\begin{thm}
Let $(X,d,\delta)$ be a strong dilatation structure. Then for any $x \in X$ the triple 
 $\displaystyle (U(x), \Sigma^{x}, \delta^{x})$ is a normed local conical group,
 with the norm induced by the distance $\displaystyle d^{x}$.
\label{tgene}
\end{thm}

The conical group $\displaystyle (U(x), \Sigma^{x}, \delta^{x})$ can be regarded as the tangent space 
of $(X,d, \delta)$ at $x$. Further will be denoted by: 
$\displaystyle T_{x} X =  (U(x), \Sigma^{x}, \delta^{x})$.

The dilatation structure on this conical group has dilatations defined by 
\begin{equation}
\bar{\delta}^{x, u}_{\varepsilon} y \, = \, \Sigma^{x}\left( u, \delta^{x}_{\varepsilon} \Delta^{x}( u, y) \right) \quad .
 \label{tandil}
\end{equation}

\subsection{Topological considerations}
 
 In this subsection we compare various topologies and uniformities related to a dilatation structure.  
  
 The axiom A3 implies that for any $x \in X$ the function $\displaystyle d^{x}$ is continuous, therefore 
 open sets with respect to $\displaystyle d^{x}$ are open with respect to $d$. 
 
 If $(X,d)$ is separable and $\displaystyle d^{x}$ is non degenerate then $\displaystyle (U(x), d^{x})$ 
 is also separable and the topologies of $d$ and $\displaystyle d^{x}$ are the same. 
 Therefore $\displaystyle (U(x), d^{x})$ is also locally compact (and a set is compact  with respect 
 to $d^{x}$ if and only if it is compact with respect to $d$). 
 
 If  $(X,d)$ is separable and $\displaystyle d^{x}$ is non degenerate then the uniformities induced by 
 $d$ and  $\displaystyle d^{x}$ are the same. Indeed, let 
 $\displaystyle \left\{u_{n} \mbox{ : } n \in \mathbb{N}\right\}$ 
 be a dense set in $U(x)$, with $\displaystyle x_{0}=x$. 
 We can embed $\displaystyle (U(x),  (\delta^{x}, \varepsilon))$ (see definition \ref{drelative}) isometrically in the separable Banach space 
 $\displaystyle l^{\infty}$, for any $\varepsilon \in (0,1)$, by the function 
 $$\phi_{\varepsilon}(u) = \left( \frac{1}{\varepsilon} d(\delta^{x}_{\varepsilon}u,  \delta^{x}_{\varepsilon}x_{n}) - \frac{1}{\varepsilon} d(\delta^{x}_{\varepsilon}x,  \delta^{x}_{\varepsilon}x_{n})\right)_{n}  . $$
 A reformulation of point (a) in theorem \ref{thcone} is that on compact sets $\displaystyle \phi_{\varepsilon}$ uniformly converges to the isometric embedding of $\displaystyle (U(x), d^{x})$ 
 $$\phi(u) = \left(  d^{x}(u,  x_{n}) - d^{x}(x, x_{n})\right)_{n}  . $$
Remark that the uniformity induced by $(\delta,\varepsilon)$ is the same as the uniformity 
induced by $d$, and that it is the same induced from the uniformity on 
$\displaystyle l^{\infty}$ by 
the embedding $\displaystyle \phi_{\varepsilon}$. We proved that the uniformities induced by 
 $d$ and  $\displaystyle d^{x}$ are the same.

From previous considerations we deduce the following characterization of  tangent
spaces associated to a dilatation structure. 

\begin{cor}
Let $(X,d,\delta)$ be a strong dilatation structure with group $\Gamma = (0,+\infty)$. 
Then for any $x \in X$ the local group  
 $\displaystyle (U(x), \Sigma^{x})$ is locally a simply connected Lie group 
 whose Lie algebra admits a positive graduation (a homogeneous group).
 \label{cortang}
\end{cor}

\paragraph{Proof.}
Use the facts:  $\displaystyle (U(x), \Sigma^{x})$ is a locally compact
group (from previous topological considerations)  which admits 
$\displaystyle \delta^{x}$ as a contractive automorphism 
group (from theorem \ref{tgene}). Apply then  Siebert proposition \ref{psiebert} ( which is \cite{siebert} proposition 5.4).
 \quad $\square$

\subsection{Differentiability with respect to dilatation structures}

We briefly explain the notion of differentiability associated 
to dilatation structures (section 7.2 \cite{buligadil1}). 
First we need the natural definition below. 

\begin{dfn}
 Let $(N,\delta)$ and $(M,\bar{\delta})$ be two  conical groups. A function $f:N\rightarrow M$ is a conical group morphism if $f$ is a group morphism and for any $\varepsilon>0$ and $u\in N$ we have 
 $\displaystyle f(\delta_{\varepsilon} u) = \bar{\delta}_{\varepsilon} f(u)$. 
\label{defmorph}
\end{dfn}

The definition of the derivative with respect to dilatations structures follows. 

 \begin{dfn}
 Let $(X, \delta , d)$ and $(Y, \overline{\delta} , \overline{d})$ be two 
 strong dilatation structures and $f:X \rightarrow Y$ be a continuous function. The function $f$ is differentiable in $x$ if there exists a 
 conical group morphism  $\displaystyle Q^{x}:T_{x}X\rightarrow T_{f(x)}Y$, defined on a neighbourhood of $x$ with values in  a neighbourhood  of $f(x)$ such that 
\begin{equation}
\lim_{\varepsilon \rightarrow 0} \sup \left\{  \frac{1}{\varepsilon} \overline{d} \left( f\left( \delta^{x}_{\varepsilon} u\right) ,  \overline{\delta}^{f(x)}_{\varepsilon} Q^{x} (u) \right) \mbox{ : } d(x,u) \leq \varepsilon \right\}Ê  = 0 , 
\label{edefdif}
\end{equation}
The morphism $\displaystyle Q^{x}$ is called the derivative of $f$ at $x$ and will be sometimes denoted by $Df(x)$.

The function $f$ is uniformly differentiable if it is differentiable everywhere and the limit in (\ref{edefdif}) 
is uniform in $x$ in compact sets. 
\label{defdif}
\end{dfn}

\section{Infinitesimal affine geometry of dilatation structures}

\subsection{Affine transformations}

\begin{dfn}
Let $(X,d,\delta)$ be a    dilatation structure. A transformation $A:X\rightarrow X$ is affine if it 
is Lipschitz and it commutes with dilatations in the following sense: for any $x\in X$, $u \in U(x)$ and 
$\varepsilon \in \Gamma$, $\nu(\varepsilon) < 1$, if  $A(u) \in U(A(x))$ then  
$$ A \delta^{x}_{\varepsilon} = \delta^{A(x)}_{\varepsilon} A(u) \quad  .$$
The local  group of affine transformations, denoted by $Aff(X,d,\delta)$ is formed by all invertible and 
bi-lipschitz affine transformations of $X$. 
\label{defgl}
\end{dfn}

$Aff(X,d,\delta)$ is  indeed a local group. In order to see this 
 we start from the  
remark that if $A$ is Lipschitz then there exists $C>0$ such that for all $x\in X$ and $u \in B(x,C)$ we have $A(u)\in U(A(x))$.  The inverse of $A \in Aff(X,d,\delta)$ is then affine. Same considerations apply for the composition of two affine, bi-lipschitz and invertible transformations. 

In the particular case of  
$X$ finite dimensional real, normed vector space, 
$d$ the distance given by the norm, $\Gamma = (0,+\infty)$ and dilatations 
$\displaystyle \delta_{\varepsilon}^{x} u = x + \varepsilon(u-x)$, 
an affine transformation in the sense of definition \ref{defgl} is an affine transformation of the vector 
space $X$.

\begin{prop}
Let $(X,d,\delta)$ be a   dilatation structure and $A:X\rightarrow X$ an affine transformation. Then: 
\begin{enumerate}
\item[(a)] for all $x\in X$, $u,v \in U(x)$ sufficiently close to $x$, we have: 
$$A \, \Sigma_{\varepsilon}^{x}(u,v)  = \Sigma_{\varepsilon}^{A(x)}(A(u),A(v)) \quad . $$
\item[(b)] for all $x\in X$, $u \in U(x)$ sufficiently close to $x$, we have: 
$$ A \, inv^{x}(u) = \, inv^{A(x)} A(u) \quad . $$
\end{enumerate}
\label{plinear}
\end{prop}

\begin{prop}
Let $(X,d,\delta)$ be a strong  dilatation structure and $A:X\rightarrow X$ an affine 
transformation. Then  $A$ is uniformly differentiable and the derivative equals $A$.
\end{prop}

The proofs are straightforward, just use the commutation with dilatations.

\subsection{Infinitesimal linearity of dilatation structures}

We begin by an explanation of the term ''sufficiently closed``, which will be used repeatedly 
in the following. 

We work in a dilatation structure $(X,d, \delta)$. Let $K \subset X$ be a compact, non empty set. Then there
is a constant $C(K) > 0$, depending on the set $K$ such that for any $\varepsilon,\mu \in 
\Gamma$ with $\nu(\varepsilon),\nu(\mu) \in (0,1]$ and any $x,y,z \in K$ with 
$d(x,y), d(x,z), d(y,z) \leq C(K)$  we have 
$$\delta^{y}_{\mu}z \in V_{\varepsilon}(x) \quad , \quad
\delta_{\varepsilon}^{x} z \in V_{\mu}(\delta^{x}_{\varepsilon} y) \quad .$$
Indeed, this is coming from the uniform (with respect to K) estimates:  
$$d(\delta^{x}_{\varepsilon} y, \delta^{x}_{\varepsilon} z) \leq \varepsilon 
d^{x}(y,z) + \varepsilon \mathcal{O}(\varepsilon) \quad , $$
$$d(x, \delta^{y}_{\mu} z) \leq d(x,y) + d(y, \delta^{y}_{\mu}z) 
\leq d(x,y) + \mu d^{y}(y,z) + \mu \mathcal{O}(\mu) \quad . $$

\begin{dfn} 
A property $\displaystyle \mathcal{P}(x_{1},x_{2},
x_{3}, ...)$ holds for $\displaystyle x_{1}, x_{2}, x_{3},
...$ sufficiently closed if for any compact, non empty set $K \subset X$, there
is a positive constant $C(K)> 0$ such that $\displaystyle \mathcal{P}(x_{1},x_{2},
x_{3}, ...)$ is true for any $\displaystyle x_{1},x_{2},
x_{3}, ... \in K$ with $\displaystyle d(x_{i}, x_{j}) \leq C(K)$.  
\end{dfn}

For example, we may say that the expressions 
$$\delta_{\varepsilon}^{x} \delta^{y}_{\mu} z \quad , \quad
\delta^{\delta^{x}_{\varepsilon} y}_{\mu} \delta^{x}_{\varepsilon} z$$ 
are well defined for any $x,y,z \in X$ sufficiently closed and for any 
$\varepsilon,\mu \in 
\Gamma$ with $\nu(\varepsilon),\nu(\mu) \in (0,1]$.

\begin{dfn}
A   dilatation structure $(X,d,\delta)$ is linear if for any $\varepsilon,\mu \in 
\Gamma$ such that $\nu(\varepsilon),\nu(\mu) \in (0,1]$, and for any 
$x,y,z \in X$ sufficiently closed we have 
$$\delta_{\varepsilon}^{x} \,  \delta^{y}_{\mu} z \ = \ 
\delta^{\delta^{x}_{\varepsilon} y}_{\mu} \delta^{x}_{\varepsilon} z \quad .$$
\label{defilin}
\end{dfn}

This definition means simply that a linear dilatation structure is a dilatation structure with the property 
that dilatations are affine in the sense of definition \ref{defgl}.

Let us  look  to a dilatation structure in finer details. We do this by defining induced dilatation structures from a given one.

\begin{dfn}
Let $(X,\delta,d)$ be a   dilatation structure and $x\in X$ a point. 
In a neighbourhood $U(x)$ of $x$, for  any $\mu\in (0,1)$ 
we define the distances:
$$(\delta^{x},\mu)(u,v) = \frac{1}{\mu} d(\delta^{x}_{\mu} u , \delta^{x}_{\mu} v) . $$
\label{drelative}
\end{dfn}

The next theorem shows  that on a dilatation structure we almost have 
translations 
(the operators $\displaystyle \Sigma^{x}_{\varepsilon}(u, \cdot)$), which are 
almost isometries (that is, not with respect to the distance $d$, but with respect to distances of type $\displaystyle (\delta^{x},\mu)$). 
It is almost as if we were working with  a normed 
conical group, only that we have to use families of distances and to make 
small shifts in the tangent space, as it is done  at the end of  the proof 
of theorem \ref{pshift}. 

\begin{thm}
Let $(X,\delta,d)$ be a   (strong) dilatation structure. 
For any $u\in U(x)$ and $v$ close to $u$  let us  define 
$$\hat{\delta}_{\mu, \varepsilon}^{x, u} v = \Sigma^{x}_{\mu}(u, \delta^{\delta^{x}_{\mu}u}_{\varepsilon} 
\Delta^{x}_{\mu}(u,v)) = \delta^{x}_{\mu^{-1}} \delta_{\varepsilon}^{\delta_{\mu}^{x} u} \delta_{\mu}^{x} v \quad  . $$
Then $\displaystyle (U(x),\hat{\delta}^{x}_{\mu}, (\delta^{x},\mu))$ is a  (strong) dilatation structure. 

The transformation 
$\displaystyle \Sigma^{x}_{\mu}(u, \cdot)$ is an isometry from 
$\displaystyle (\delta^{\delta^{x}_{\mu} u}, \mu)$ to 
$\displaystyle (\delta^{x}, \mu)$. Moreover, we have 
$\displaystyle  \Sigma^{x}_{\mu} (u, \delta_{\mu}^{x} u) = u$.
\label{pshift}
\end{thm}

\paragraph{Proof.}
We have to check the axioms.  The first part of axiom A0 is an easy consequence of theorem 
\ref{thcone} for $(X,\delta,d)$. The second part of A0, A1 and A2 are true based on simple computations. 

The first interesting fact is related to axiom A3. Let us compute, for $v,w \in U(x)$, 
$$\frac{1}{\varepsilon} (\delta^{x},\mu)(\hat{\delta}^{x, u}_{\mu, \varepsilon} v, \hat{\delta}^{x, u}_{\mu \varepsilon} w) = 
\frac{1}{\varepsilon \mu} d( \delta^{x}_{\mu} \hat{\delta}^{x, u}_{\mu \varepsilon} v, \delta^{x}_{\mu}
\hat{\delta}^{x, u}_{\mu \varepsilon} w) = $$
$$ = \frac{1}{\varepsilon \mu} d( \delta_{\varepsilon}^{\delta^{x}_{\mu} u} \delta_{\mu}^{x} v ,  
\delta_{\varepsilon}^{\delta^{x}_{\mu} u} \delta_{\mu}^{x} w) = \frac{1}{\varepsilon\mu} d(  
\delta_{\varepsilon \mu}^{\delta^{x}_{\mu} u} \Delta_{\mu}^{x}(u,v),  \delta_{\varepsilon \mu}^{\delta^{x}_{\mu} u} \Delta_{\mu}^{x}(u,w)) =  $$
$$= (\delta^{\delta_{\mu}^{x} u}, \varepsilon \mu) (  \Delta_{\mu}^{x}(u,v) ,  \Delta_{\mu}^{x}(u,w)) . $$
The axiom A3 is then a consequence of axiom A3 for $(X,d, \delta)$ and we have 
$$\lim_{\varepsilon\rightarrow 0} \frac{1}{\varepsilon} (\delta^{x},\mu)(\hat{\delta}^{x, u}_{\mu \varepsilon} v, \hat{\delta}^{x, u}_{\mu \varepsilon} w) = d^{\delta_{\mu}^{x} u} (  \Delta_{\mu}^{x}(u,v) ,  \Delta_{\mu}^{x}(u,w)) . $$
The axiom A4 is also a straightforward consequence of A4 
for $(X,d, \delta)$. The second part of the theorem is a simple computation.  \quad $\square$

The induced dilatation structures $\displaystyle (U(x),\hat{\delta}^{x}_{\mu}, (\delta^{x},\mu))$ should converge in some sense to the dilatation structure 
on the tangent space at $x$, as $\nu(\mu)$ converges to zero. Remark that we have one easy convergence in strong dilatation structures: 
$$\lim_{\mu \rightarrow 0}  \hat{\delta}^{x, u}_{\mu, \varepsilon} v \, = \, \bar{\delta}^{x, u}_{\varepsilon} v $$ 
where $\displaystyle \bar{\delta}^{x}$ are the dilatations in the tangent space at $x$, cf. (\ref{tandil}). 
Indeed, this comes from: 
$$ \hat{\delta}^{x, u}_{\mu, \varepsilon} v \, = \, \Sigma^{x}_{\mu}(u, \delta^{\delta^{x}_{\mu}u}_{\varepsilon} 
\Delta^{x}_{\mu}(u,v)) $$ 
so, when $\nu(\mu)$ converges we get the mentioned limit. 

The following proposition gives  a more precise estimate: the order of approximation of the dilatations $\delta$ by dilatations 
$\displaystyle \hat{\delta}^{x}_{\varepsilon}$, in neighbourhoods of $\displaystyle \delta^{x}_{\varepsilon} y$ of order $\varepsilon$, 
as $\nu(\varepsilon)$ goes to zero.

\begin{prop}
Let $(X,\delta,d)$ be a   dilatation structure. With the notations of theorem \ref{pshift} we introduce 
$$\hat{\delta}_{\varepsilon}^{x, u} v =  \hat{\delta}_{\varepsilon, \varepsilon}^{x, u} v  = \delta^{x}_{\varepsilon^{-1}} \delta_{\varepsilon}^{\delta_{\varepsilon}^{x} u} \delta_{\varepsilon}^{x} v \quad  . $$
Then we have for any $x, y, v$ sufficiently closed: 
\begin{equation}
 \lim_{\varepsilon \rightarrow 0} \frac{1}{\varepsilon} \, (\delta^{x}, \varepsilon) \left(  \delta_{\varepsilon}^{\delta_{\varepsilon}^{x} y} v \, , \, 
\hat{\delta}_{\varepsilon}^{x, \delta_{\varepsilon}^{x} y} v \right) \, = \, 0 \quad . 
\label{firlin}
\end{equation}
\label{plin1}
\end{prop}

\paragraph{Proof.}
We start by a computation: 
$$ \frac{1}{\varepsilon} \, (\delta^{x}, \varepsilon) \left(  \delta_{\varepsilon}^{\delta_{\varepsilon}^{x} y} v \, , \, 
\hat{\delta}_{\varepsilon}^{x,\delta_{\varepsilon}^{x} y} v \right) = \frac{1}{\varepsilon^{2}} \, d \, \left( \delta^{x}_{\varepsilon} \delta_{\varepsilon}^{\delta_{\varepsilon}^{x} y} v \, , \, \delta^{x}_{\varepsilon} 
\hat{\delta}_{\varepsilon}^{x, \delta_{\varepsilon}^{x} y} v \right) = $$ 
$$ = \frac{1}{\varepsilon^{2}} \, d \, \left( \delta^{x}_{\varepsilon^{2}} \Sigma_{\varepsilon}^{x}(y, v) \, , \, \delta^{x}_{\varepsilon^{2}} \,  \delta^{x}_{\varepsilon^{-2}} \delta_{\varepsilon^{2}}^{\delta_{\varepsilon^{2}}^{x} y} \, \Delta_{\varepsilon}^{x} (\delta^{x}_{\varepsilon} y , v) \right) = $$  
$$ =  \frac{1}{\varepsilon^{2}} \, d \, \left( \delta^{x}_{\varepsilon^{2}} \Sigma_{\varepsilon}^{x}(y, v) \, , \, \delta^{x}_{\varepsilon^{2}} 
\Sigma^{x}_{\varepsilon^{2}}\left( y, \Delta_{\varepsilon}^{x} (\delta^{x} y , v) \right)  \right) \quad . $$ 
This last expression converges as $\nu(\varepsilon)$ goes to $0$ to 
$$d^{x} \left( \Sigma^{x}(y,v) , \Sigma^{x}(y , \Delta^{x}(x,v)) \right) = d^{x} \left( v, \Delta^{x}(x,v) \right) = 0 $$
\quad $\square$

The result from this proposition  indicates that strong dilatation structures are infinitesimally linear. In order to make a precise statement 
we need a measure for nonlinearity 
of a dilatation structure, given in the next definition. Then we have to repeat the computations from the proof of proposition \ref{plin1} in a slightly 
different setting, related to this measure of nonlinearity.

\begin{dfn}
The following expression: 
\begin{equation}
Lin(x,y,z; \varepsilon, \mu) \ = \ d\left( \delta_{\varepsilon}^{x} \, 
\delta^{y}_{\mu}z \, , \,  \delta^{\delta^{x}_{\varepsilon} y}_{\mu}
\delta^{x}_{\varepsilon} z\right) 
\label{deflinfunct}
\end{equation}
is a measure of lack of linearity, for a general    dilatation structure. 
\end{dfn}

The next theorem shows that indeed, infinitesimally any  strong dilatation
structure is linear.

\begin{thm}
Let $(X,d,\delta)$ be a  strong dilatation structure. Then for any $x, y, z \in X$ 
sufficiently close we have 
\begin{equation}
\lim_{\varepsilon \rightarrow 0} \frac{1}{\varepsilon^{2}} \, 
Lin(x,\delta_{\varepsilon}^{x} y,  z ; \varepsilon,
\varepsilon) \ = \ 0 \quad . 
\label{inflin}
\end{equation}
\label{thinflin}
\end{thm}

\paragraph{Proof.}
From the hypothesis of the theorem we have:
$$\frac{1}{\varepsilon^{2}} \, 
Lin(x,\delta_{\varepsilon}^{x} y, \delta_{\varepsilon}^{x} z ; \varepsilon,
\varepsilon) \ = \ \frac{1}{\varepsilon^{2}} \, d\left( 
\delta^{x}_{\varepsilon} \, \delta^{\delta^{x}_{\varepsilon} y}_{\varepsilon} z 
\, , \, 
\delta_{\varepsilon}^{\delta_{\varepsilon^{2}}^{x} y} \delta_{\varepsilon}^{x} z
\right) \ = \ $$
$$= \ \frac{1}{\varepsilon^{2}} \, d\left( \delta^{x}_{\varepsilon^{2}} \, 
\Sigma^{x}_{\varepsilon}(y,z) \, , \, \delta^{x}_{\varepsilon^{2}} \, 
\delta^{x}_{\varepsilon^{-2}} \, 
\delta_{\varepsilon}^{\delta_{\varepsilon^{2}}^{x} y} \delta_{\varepsilon}^{x} z
\right) \ = \ $$ 
$$= \ \frac{1}{\varepsilon^{2}} \, d\left( \delta^{x}_{\varepsilon^{2}}
\, \Sigma^{x}_{\varepsilon}(y,z) \, , \, \delta^{x}_{\varepsilon^{2}} \, 
\Sigma^{x}_{\varepsilon^{2}} (y \, , \, 
\Delta^{x}_{\varepsilon}(\delta^{x}_{\varepsilon} y ,  z))\right) \ = \ $$ 
$$= \ \mathcal{O}(\varepsilon^{2}) + d^{x}\left( \Sigma^{x}_{\varepsilon}(y,z) \,
, \, \Sigma^{x}_{\varepsilon^{2}} (y \, , \, 
\Delta^{x}_{\varepsilon}(\delta^{x}_{\varepsilon} y ,  z))\right) \quad .$$
The dilatation structure satisfies A4, therefore as $\varepsilon$ goes to $0$ we
obtain: 
$$\lim_{\varepsilon \rightarrow 0} \frac{1}{\varepsilon^{2}} \, 
Lin(x,\delta_{\varepsilon}^{x} y, \delta_{\varepsilon}^{x} z ; \varepsilon,
\varepsilon) \ = \ d^{x}\left( \Sigma^{x}(y,z) \,
, \, \Sigma^{x} (y \, , \, 
\Delta^{x}(x ,  z))\right) \ =$$ 
$$= \  d^{x}\left( \Sigma^{x}(y,z) \,
, \, \Sigma^{x}(y , z)\right) \ =  \ 0 \quad . \quad \quad \quad \square$$

\subsection{Linear  strong dilatation structures}

Remark that for general dilatation structures the "translations" 
$\displaystyle \Delta^{x}_{\varepsilon}(u, \cdot)$  are not affine. 
Nevertheless, they commute with dilatation in a known way: for any $u,v$ sufficiently close to $x$ and $\mu \in \Gamma$, 
$\nu(\mu) < 1$,     we have: 
$$\Delta^{x}_{\varepsilon} \left( \delta^{x}_{\mu} u,  \delta^{x}_{\mu} v \right) =  \delta^{ \delta^{x}_{\epsilon \mu} u}_{\mu} \Delta^{x}_{\varepsilon \mu} (u,v) \quad . $$

 This is important, because the 
transformations $\displaystyle \Sigma^{x}_{\varepsilon}(u, \cdot)$ really behave as translations.  
The reason for which such transformations are not affine is that dilatations are
generally not affine.

Linear dilatation structures are very particular dilatation structures. The next
proposition gives a family of examples of linear dilatation structures. 

\begin{prop}
The dilatation structure associated to a normed conical group is linear.
\label{pexlin}
\end{prop}

\paragraph{Proof.}
Indeed, for the dilatation structure associated to a normed conical group we
have, with the notations from definition \ref{defilin}: 
$$\delta^{\delta^{x}_{\varepsilon} y}_{\mu}
\delta^{x}_{\varepsilon} z \ = \ \left(x \delta_{\varepsilon} (x^{-1}y) \right) 
\, \delta_{\mu} \left( \delta_{\varepsilon} (y^{-1}x) \, x^{-1} \, x \, 
\delta_{\varepsilon} (x^{-1} z) \right) \ =  $$
$$= \  \left(x \delta_{\varepsilon} (x^{-1}y) \right) 
\, \delta_{\mu} \left( \delta_{\varepsilon} (y^{-1}x)  \, 
\delta_{\varepsilon} (x^{-1} z) \right) \ = \ 
\left(x \delta_{\varepsilon} (x^{-1}y) \right) 
\, \delta_{\mu} \left( \delta_{\varepsilon} (y^{-1}z) \right) \ =  $$
$$= \ x \left(\delta_{\varepsilon} (x^{-1}y) 
\, \delta_{\varepsilon} \, \delta_{\mu} (y^{-1}z) \right) \ = \ 
x \, \delta_{\varepsilon} \left(x^{-1}y 
\,  \delta_{\mu} (y^{-1}z)\right)  \ = \ \delta_{\varepsilon}^{x} \, 
\delta^{y}_{\mu}z \quad .$$
Therefore the dilatation structure is linear. \quad $\square$

The affinity of translations  
$\displaystyle \Sigma^{x}_{\varepsilon}$  is related to the linearity of the
dilatation structure, as described in the theorem below, point (a).  
As a consequence, we prove at point (b) that a linear and strong
dilatation structure comes from a conical group.

\begin{thm}
Let $(X,d,\delta)$ be a    dilatation structure. 
\begin{enumerate}
\item[(a)] If the dilatation structure is linear  then all    transformations 
$\displaystyle \Delta^{x}_{\varepsilon}(u, \cdot)$ are affine for any $u \in X$. 
\item[(b)]   If the dilatation structure is strong (satisfies A4) then it is 
linear  if and only if the dilatations  come  from the  dilatation structure of a 
 conical group, precisely for any $x \in X$ there is an open neighbourhood $D \subset X$
 of $x$ such that  $\displaystyle (\overline{D}, d^{x}, \delta)$ is the same 
 dilatation structure as the dilatation structure of the tangent space 
 of $(X,d,\delta)$ at $x$.
\end{enumerate}
\label{tdilatlin}
\end{thm}

\paragraph{Proof.}
(a) If dilatations are affine, then let $\varepsilon, \mu \in \Gamma$, $\nu(\varepsilon), 
\nu(\mu) \leq 1$, and $x, y, u, v \in X$  such that the following computations make sense. We have: 
$$\Delta^{x}_{\varepsilon} (u, \delta_{\mu}^{y} v ) = \delta_{\varepsilon^{-1}}^{\delta^{x}_{\varepsilon} u} \delta_{\varepsilon}^{x} \delta_{\mu}^{y} v \quad . $$
Let $\displaystyle A_{\varepsilon} = \delta_{\varepsilon^{-1}}^{\delta^{x}_{\varepsilon} u}$. We compute: 
$$\delta_{\mu}^{\Delta_{\varepsilon}^{x}(u,y)} \Delta_{\varepsilon}^{x} (u,v)  = \delta_{\mu}^{A_{\varepsilon} \delta_{\varepsilon}^{x} y} A_{\varepsilon} \delta_{\varepsilon}^{x} v \quad . $$
We use twice the affinity of dilatations:  
$$ \delta_{\mu}^{\Delta_{\varepsilon}^{x}(u,y)} \Delta_{\varepsilon}^{x} (u,v)  = A_{\varepsilon} 
\delta_{\mu}^{\delta_{\varepsilon}^{x} y} \delta_{\varepsilon}^{x} v = 
 \delta_{\varepsilon^{-1}}^{\delta^{x}_{\varepsilon} u} \delta_{\varepsilon}^{x} \delta_{\mu}^{y} v \quad . $$
We proved that: 
$$\Delta^{x}_{\varepsilon} (u, \delta_{\mu}^{y} v ) =  \delta_{\mu}^{\Delta_{\varepsilon}^{x}(u,y)} \Delta_{\varepsilon}^{x} (u,v)  \quad , $$
which is the conclusion of   the part (a).

(b) Suppose that the dilatation structure is strong. If dilatations are affine, then by point (a) the  transformations $\displaystyle
\Delta^{x}_{\varepsilon}(u, \cdot)$ are affine as well for any $u \in X$. Then, with notations made before, for $y = u$ we get 
$$\Delta^{x}_{\varepsilon} (u, \delta_{\mu}^{u} v ) =  \delta_{\mu}^{\delta_{\varepsilon}^{x} u} \Delta_{\varepsilon}^{x} (u,v)  \quad , $$
which implies 
$$\delta_{\mu}^{u} v = \Sigma^{x}_{\varepsilon} ( u, \delta^{x}_{\mu} \Delta_{\varepsilon}^{x}(u,v)) \quad . $$
We pass to the limit with $\varepsilon \rightarrow 0$ and we obtain: 
$$ \delta_{\mu}^{u} v = \Sigma^{x}( u, \delta^{x}_{\mu} \Delta^{x}(u,v)) \quad . $$
We recognize at the  right hand side the dilatations associated to the conical group 
$\displaystyle T_{x} X$. 

By proposition \ref{pexlin} the opposite implication is straightforward, 
because the dilatation structure of any conical group is linear. 
\quad $\square$

\section{Noncommutative affine geometry}

We propose here to call ''noncommutative affine geometry`` the  generalization of affine geometry described in theorem 
\ref{taffine}, but without the restriction $\Gamma = (0,+\infty)$. For short, noncommutative affine geometry in the sense explained further 
is the study of the properties of linear strong dilatation structures. Equally, by theorem \ref{tdilatlin}, it is the study of normed conical 
groups.

As a motivation for this name, in the proposition below we give a relation, true for linear dilatation
structures, with an interesting  interpretation. We shall explain what this relation means 
in the most trivial case: the dilatation structure associated to a real normed  affine space. In this case, 
 for any  points $x,u,v$,  let us denote 
$\displaystyle w = \Sigma^{x}_{\varepsilon}(u,v)$. Then $w$ equals (approximatively, due to the parameter 
$\varepsilon$)  the sum   $\displaystyle u+_{x} v$.  
Denote also  $\displaystyle w' = 
\Delta^{u}_{\varepsilon}(x,v)$; then  $w'$ is 
 (approximatively) equal to  the difference 
between  $v$ and $x$ based at $u$. In our space (a classical  affine space 
over a vector space)  we  have $\displaystyle w = w'$. The next proposition shows that 
the  same is true for any linear dilatation
structure. 

\begin{prop}
For a linear    dilatation structure $(X,\delta,d)$, for any $x,u,v \in X$ 
sufficiently closed and for any $\varepsilon \in \Gamma$, $\nu(\varepsilon) \leq
1$, we have: 
$$\Sigma^{x}_{\varepsilon}(u,v) \ = \ \Delta^{u}_{\varepsilon}(x,v) \quad . $$
\end{prop}

\paragraph{Proof.}
We have the following string of equalities, by using twice the linearity of the dilatation
structure: 
$$\Sigma_{\varepsilon}^{x}(u,v) \ = \ \delta_{\varepsilon^{-1}}^{x}
\delta_{\varepsilon}^{\delta_{\varepsilon}^{x} u} v \ = \
\delta_{\varepsilon}^{u} \, \delta_{\varepsilon^{-1}}^{x} v \ = \ $$
$$= \ \delta_{\varepsilon^{-1}}^{\delta_{\varepsilon}^{u} x}
\delta_{\varepsilon}^{u} v \ = \ \Delta_{\varepsilon}^{u}(x,v) \quad . $$
The proof is done. \quad $\square$

\subsection{Inverse semigroups and Menelaos theorem}

Here we prove that for strong  dilatation structures linearity is equivalent 
to  a generalization of the statement from corollary \ref{corunu}. 
The result is new  for Carnot groups and the proof seems to be new even for 
vector spaces.

\begin{dfn}
A semigroup $S$ is an inverse semigroup if for any $x \in S$ there is 
an unique element $\displaystyle x^{-1} \in S$ such that 
$\displaystyle x \, x^{-1}  x = x$ and $\displaystyle x^{-1} x \,  x^{-1}  =
x^{-1}$. 
\label{dinvsemi}
\end{dfn}

An important  example of an inverse semigroup is $I(X)$, the class of all
bijective maps $\phi: \, dom \, \phi \, \rightarrow \, im \, \phi$, with 
$dom \, \phi , \, im \, \phi \subset X$. The semigroup operation is the
composition of functions in the largest domain where this makes sense. 

By the Vagner-Preston representation theorem \cite{howie} every inverse
semigroup is isomorphic to a subsemigroup of $I(X)$, for some set $X$.

\begin{dfn}
A dilatation structure $(X,d,\delta)$ has the Menelaos property if 
for any two sufficiently closed $x,y \in X$  and for any 
$\varepsilon,\mu \in \Gamma$ with $\nu(\varepsilon),\nu(\mu) \in (0,1)$ we have 
$$\delta^{x}_{\varepsilon} \, \delta^{y}_{\mu} \ = \ \delta^{w}_{\varepsilon \mu}
\quad , $$
where $w \in X$ is the fixed point of the contraction $\displaystyle 
\delta^{x}_{\varepsilon} \delta^{y}_{\mu}$ (thus depending on $x,y$ and 
$\varepsilon , \mu$). 
\label{defmene}
\end{dfn}

\begin{thm}
A linear dilatation structure has the Menelaos property. 
\label{thmenelin}
\end{thm}

\paragraph{Proof.} Let $x, y \in X$ be sufficiently closed and 
$\varepsilon,\mu \in \Gamma$ with $\nu(\varepsilon),\nu(\mu) \in (0,1)$. 
We shall define two sequences $\displaystyle x_{n}, y_{n} \in X$, $n \in
\mathbb{N}$. 

We begin with $\displaystyle x_{0} = x$, $y_{0} = y$. Suppose further that 
$\displaystyle x_{n}, y_{n}$ were defined and that they are sufficiently closed.
Then we use twice the linearity of the dilatation structure: 
$$\delta_{\varepsilon}^{x_{n}} \, \delta^{y_{n}}_{\mu} \ = \ 
\delta_{\mu}^{\delta_{\varepsilon}^{x_{n}} y_{n}} \, 
\delta_{\varepsilon}^{x_{n}} \ = \ \delta_{\varepsilon}^{\delta_{\mu}^{
\delta_{\varepsilon}^{x_{n}} y_{n}} x_{n}} \, 
\delta_{\mu}^{\delta_{\varepsilon}^{x_{n}} y_{n}} \quad . $$
We shall define then by induction 
\begin{equation}
x_{n+1} \ = \ \delta_{\mu}^{\delta_{\varepsilon}^{x_{n}} y_{n}} x_{n} \quad , 
\quad y_{n+1} \ = \ \delta_{\varepsilon}^{x_{n}} y_{n} \quad .
\label{definduc}
\end{equation}
Provided that we prove by induction that $\displaystyle x_{n}, y_{n}$ are
sufficiently closed, we arrive at the conclusion  that for any 
$n \in \mathbb{N}$ 
\begin{equation}
\delta_{\varepsilon}^{x_{n}} \, \delta^{y_{n}}_{\mu} \ = \
\delta_{\varepsilon}^{x} \, \delta^{y}_{\mu} \quad . 
\label{equlimain}
\end{equation}
The points $\displaystyle x_{0}, y_{0}$ are sufficiently closed by hypothesis. 
Suppose now that $\displaystyle x_{n}, y_{n}$ are sufficiently closed. 
Due to the linearity of the dilatation structure, 
 we can write the first part of (\ref{definduc}) as: 
$$x_{n+1} \ = \delta_{\varepsilon}^{x_{n}} \, \delta^{y_{n}}_{\mu} x_{n} \quad . $$
Then we can estimate the distance between $\displaystyle x_{n+1}, y_{n+1}$ like this: 
$$d(x_{n+1},y_{n+1}) \ = \ d( \delta_{\varepsilon}^{x_{n}} \, \delta^{y_{n}}_{\mu} x_{n}
, \delta_{\varepsilon}^{x_{n}} y_{n} ) \ = \ \nu(\varepsilon) \, 
d(\delta^{y_{n}}_{\mu} x_{n} , y_{n}) \ = \ \nu(\varepsilon \mu) d(x_{n}, y_{n})
\quad .$$
From $\nu(\varepsilon \mu) < 1$ it follows that $\displaystyle x_{n+1}, y_{n+1}
$ are sufficiently closed. By induction we deduce that for all 
$n \in \mathbb{N}$ the points $\displaystyle x_{n+1}, y_{n+1}
$ are sufficiently closed. We also find out that 
\begin{equation}
\lim_{n \rightarrow \infty} d(x_{n}, y_{n}) \ = 
\ 0 \quad . 
\label{equlim}
\end{equation}

From relation (\ref{equlimain}) we deduce that  the first part of 
(\ref{definduc}) can be written as: 
$$x_{n+1} \ = \delta_{\varepsilon}^{x_{n}} \, \delta^{y_{n}}_{\mu} x_{n} \ = \ 
\delta_{\varepsilon}^{x} \, \delta^{y}_{\mu} x_{n} \quad . $$
The transformation $\displaystyle \delta_{\varepsilon}^{x} \, \delta^{y}_{\mu}$ 
is a contraction of coefficient  $\nu(\varepsilon \mu) < 1$, therefore we easily
get: 
\begin{equation}
\lim_{n \rightarrow \infty} x_{n} \ = \ w \quad , 
\label{equlim2}
 \end{equation}
where $w$ is the unique fixed point of the contraction   $\displaystyle
\delta_{\varepsilon}^{x} \, \delta^{y}_{\mu}$. 

We put together (\ref{equlim}) and (\ref{equlim2}) and we get the limit: 
\begin{equation}
\lim_{n \rightarrow \infty} y_{n} \ = \ w \quad , 
\label{equlim3}
 \end{equation}
Using relations (\ref{equlim2}), (\ref{equlim3}), we may pass to the limit with $n \rightarrow \infty$ 
 in relation (\ref{equlimain}): 
 $$ \delta_{\varepsilon}^{x} \, \delta^{y}_{\mu} \ = \ \lim_{n \rightarrow 
 \infty} \delta_{\varepsilon}^{x_{n}} \, \delta^{y_{n}}_{\mu} \ = \ 
 \delta_{\varepsilon}^{w} \, \delta^{w}_{\mu} \ = \ \delta^{w}_{\varepsilon 
 \mu} \quad .$$ 
 The proof is done. \quad $\square$

\begin{cor}
Let $(X,d,\delta)$ be a strong linear dilatation structure, 
with group $\Gamma$ and the morphism $\nu$ injective. Then any element of the  
inverse subsemigroup of $I(X)$ generated by  dilatations is locally 
a dilatation $\displaystyle \delta^{x}_{\varepsilon}$  or  a  left translation 
 $\displaystyle \Sigma^{x}(y, \cdot)$. 
\label{cordoi}
\end{cor}

\paragraph{Proof.}
Let $(X,d,\delta)$ be a strong linear dilatation structure. From the linearity 
and theorem \ref{thmenelin} we deduce that we have to care only about
the results of compositions of two dilatations which are isometries. 
 
The dilatation structure is strong, therefore  by theorem \ref{tdilatlin} the
dilatation structure is locally coming from a conical group. 

Let us compute a composition of dilatations $\displaystyle \delta^{x}_{\varepsilon}
\delta^{y}_{\mu}$, with $\nu(\varepsilon \mu) = 1$. Because the morphism $\nu$ is injective, 
it follows that $\displaystyle \mu = \varepsilon^{-1}$. 
In a conical group 
we can make the  following computation (here $\displaystyle \delta_{\varepsilon} =
\delta^{e}_{\varepsilon}$ with $e$ the neutral element of the conical group): 
$$ \delta^{x}_{\varepsilon} \delta^{y}_{\varepsilon^{-1}} z \ = \ 
x \delta_{\varepsilon} \left( x^{-1} y \delta_{\varepsilon^{-1}} \left(y^{-1} z 
\right)\right) \ = \ x \delta_{\varepsilon} \left( x^{-1} y \right) y^{-1} z
\quad . $$
Therefore the composition of dilatations $\displaystyle \delta^{x}_{\varepsilon}
\delta^{y}_{\mu}$, with $\varepsilon \mu = 1$, is a left translation. 

Another easy computation shows that composition of left translations with
dilatations are dilatations. The proof end by remarking that all the statements 
are local. \quad $\square$

\paragraph{A counterexample.} 
\label{counter1}
The Corollary \ref{cordoi} is not true without the injectivity assumption on 
$\nu$. Indeed, consider the Carnot group $N = \mathbb{C} \times \mathbb{R}$ with the 
elements denoted by $X \in N$, $X = (x,x')$, with $x \in \mathbb{C}$, $x' \in \mathbb{R}$, and operation
$$X \, Y \ = \ (x,x') (y, y') \ = \ (x + y, x' + y' + \frac{1}{2} \,Im \, x \bar{y} )$$
We take  $\Gamma = \mathbb{C}^{*}$ and morphism $\nu: \Gamma \rightarrow (0, + \infty)$, $\nu(\varepsilon) = \mid \varepsilon \mid$. 
Dilatations are defined as: for any $\varepsilon \in \mathbb{C}^{*}$ and $X = (x,x') \in N$: 
$$\delta_{\varepsilon} X \ =  \ (\varepsilon x , \mid \varepsilon \mid^{2} x')$$
These dilatations induce the field of dilatations $\displaystyle \delta^{X}_{\varepsilon} Y \, = \, X \delta_{\varepsilon} (X^{-1} Y)$. 

The morphism $\nu$ is not injective. Let now $\varepsilon, \mu \in \mathbb{C}^{*}$ with $\varepsilon \mu = -1$ and 
$\varepsilon \in (0,1)$. An elementary (but a bit long) computation shows that for $X = (0,0)$ and $Y = (y, y')$ with 
$y \not = 0$, $y' \not = 0$,   the composition of dilatations 
$\displaystyle \delta^{X}_{\varepsilon} \delta^{Y}_{\mu}$ is not a left translation in the group $N$, nor a dilatation. \quad $\square$

Further we shall suppose that the morphism $\nu$ is always injective, if not explicitly stated otherwise. Therefore we shall 
consider $\Gamma \subset (0,+\infty)$ as a subgroup. 

\subsection{On the barycentric condition}

The barycentric condition is (Af3):  for any $\varepsilon \in (0,1)$ 
 $\displaystyle \delta_{\varepsilon}^{x} \, y \, = \delta_{1-\varepsilon}^{y} \, x$. 
In particular, the condition (Af3) tells that the transformation $\displaystyle y \mapsto \delta_{\varepsilon}^{y} x$ is also a dilatation. Is this true 
for linear dilatation structures? Theorem \ref{taffine} indicates that (Af3) is true if and only if this is a dilatation structure of a 
normed real affine space. 

\begin{prop}
Let $X$  be a normed conical group with neutral element $e$, dilatations $\delta$ and distance $d$ induced by the 
homogeneous norm $\| \cdot \|$, and $\varepsilon \in (0,1) \cap \Gamma $. Then the function 
$$ h_{\varepsilon}: X \rightarrow X \quad , \quad  \displaystyle h_{\varepsilon} (x) = x \delta_{\varepsilon} (x^{-1}) = \delta^{x}_{\varepsilon} \, e$$ is invertible and the inverse $\displaystyle g_{\varepsilon}$  has the expression  
$$g_{\varepsilon} (y) \ = \ \prod_{k = 0}^{\infty} \delta_{\varepsilon^{k}} (y) \ = \ \lim_{N \rightarrow \infty} \prod_{k = 0}^{N} \delta_{\varepsilon^{k}} (y)$$
\label{dilinv}
\end{prop}

\begin{rmk}
As the choice of the neutral element is not important,  the previous proposition says that for any $\varepsilon \in (0,1)$ and any fixed $y \in X$ 
the function $x \mapsto \delta^{x}_{\varepsilon} y$ is invertible. 
\end{rmk}

\paragraph{Proof.}
Let $\varepsilon \in (0,1)$ be fixed. For any  natural number $N$ we define $\displaystyle g_{N}: X \rightarrow X$ by 
$$g_{N}(y) \ = \ \prod_{k = 0}^{N} \delta_{\varepsilon^{k}} (y)$$
For fixed $y \in X$ $\displaystyle (g_{N}(y))_{N}$ is a Cauchy sequence. Indeed, for any $N \in \mathbb{N}$ we have: 
$$d(g_{N}(y) , g_{N+1}(y)) \ =  \ \| \delta_{\varepsilon^{N+1}} (y) \|$$
thus for any $N, M \in \mathbb{N}$, $M \geq 1$ we have 
$$d(g_{N}(y), g_{N+M}(y)) \ \leq \ \left( \sum_{k = N+1}^{M} \varepsilon^{k} \right) \, \| y \| \ \leq \ \frac{\varepsilon^{N+1}}{1-\varepsilon} \, \| y \|$$
Let then $\displaystyle g_{\varepsilon}(y) \ = \ \lim_{N \rightarrow \infty} g_{N} (y)$. We prove that $\displaystyle g_{\varepsilon}$ is the 
inverse of $\displaystyle h_{\varepsilon}$. We have, for any natural number $N$ and $y \in X$ 
$$y \, \delta_{\varepsilon} g_{N}(y) \ = \ g_{N+1}(y)$$
By passing to the limit with $N$ we get that $\displaystyle h_{\varepsilon} \circ g_{\varepsilon} (y) = y$ for any $ y \in X$. 

Let us now compute 
$$\displaystyle   g_{N} \circ  h_{\varepsilon} (x)\ = \ 
\prod_{k = 0}^{N} \delta_{\varepsilon^{k}} (x \delta_{\varepsilon} (x^{-1})) \  
= \  \prod_{k = 0}^{N} \delta_{\varepsilon^{k}} (x) \, \delta_{\varepsilon^{k+1}} (x^{-1}) 
\ = $$ 
$$ = \ x \, \delta_{\varepsilon^{N+1}} (x^{-1})$$
 therefore as $N$ goes to 
infinity we get $\displaystyle g_{\varepsilon} \circ h_{\varepsilon} (x) \ = \ x$. \quad $\square$

For any $\varepsilon \in (0,1)$ the functions $\displaystyle h_{\varepsilon}, g_{\varepsilon}$ are homogeneous, that is 
$$h_{\varepsilon} (\delta_{\mu} x) \ = \ \delta_{\mu} \, h_{\varepsilon}(x) \quad , \quad g_{\varepsilon} (\delta_{\mu} y) \ = \ \delta_{\mu} \, g_{\varepsilon}(y) $$
for any $\mu > 0$ and $x, y \in X$. 

In the presence of the barycentric condition we get the following: 

\begin{cor}
Let $(X,d,\delta)$ be a strong dilatation structure with group $\Gamma \subset (0,+\infty)$, which satisfies the barycentric condition (Af3). Then for any 
$u, v \in X$ and $\varepsilon \in (0,1) \cap \Gamma$ the points $\displaystyle inv^{u}(v)$, $u$ and $\delta^{u}_{\varepsilon}v$ are collinear in the sense:
$$d(inv^{u}(v) , u) + d(u , \delta^{u}_{\varepsilon}v ) \ = \ d(inv^{u}(v) , \delta^{u}_{\varepsilon}v )$$
\label{corbari}
\end{cor}

\paragraph{Proof.}
There is no restriction to  work with the group operation with neutral element $e$ and denote $\displaystyle \delta_{\varepsilon} := \delta_{\varepsilon}^{e}$. 
 With the 
notation from the proof of the proposition \ref{dilinv}, we use the expression 
of the function $\displaystyle g_{\varepsilon}$, we apply the homogeneous 
norm $\| \cdot \|$ and we obtain: 
$$\| g_{\varepsilon}(y) \| \ \leq \ \left( \sum_{k = 0}^{\infty} \varepsilon^{k} \right) \| x \| \ = \ \frac{1}{1-\varepsilon} \, \|y\|$$
with equality if and only if $e$, $y$ and $\displaystyle y \delta_{\varepsilon} y$ are collinear in the sense 
$\displaystyle d(e, y) + d(y, y \delta_{\varepsilon} y) \ = \ d(e, y \delta_{\varepsilon} y)$.  
The barycentric condition can be written as: $\displaystyle h_{\varepsilon}(x) = \delta_{1-\varepsilon}(x)$. We have therefore: 
$$\| x \| \ =  \ \| g_{\varepsilon} \circ h_{\varepsilon}(x) \| \  \leq \ \frac{1}{1-\varepsilon} \| h_{\varepsilon}(x)\| \ = \ 
\frac{1-\varepsilon}{1-\varepsilon} \, \|x\| \ = \ \|x \|$$ 
therefore $e$, $x$ and $\displaystyle x \delta_{\varepsilon} x$ are on a geodesic. This is true also for the choice: $\displaystyle e = inv^{u}(v)$, 
$x = u$, which gives the conclusion. \quad $\square$

We can actually say more in the case $\Gamma = (0, + \infty)$. 

\begin{prop}
 Let $(X,d,\delta)$ be a strong dilatation structure with group $\Gamma = (0,+\infty)$, which satisfies the barycentric condition (Af3). Then for 
any $x \in X$ the group operation $\displaystyle \Sigma^{x}$ is abelian and moreover the graduation of $X$, as a homogeneous group with respect to 
the operation $\displaystyle \Sigma^{x}$ has only one level.
\label{probari}
\end{prop}

\paragraph{Proof.}
Let us denote the  neutral element by  $e$ instead of $x$  and denote $\displaystyle \delta_{\varepsilon} := \delta_{\varepsilon}^{e}$. 
According to corollary \ref{cortang} $X$ is a Lie homogeneous group. The barycentric condition implies: for any $x \in X$ and $\varepsilon \in (0,1)$ 
 we have $\displaystyle \delta_{1-\varepsilon} y \ = \ y \delta_{\varepsilon} y^{-1}$, which implies: 
$$\delta_{1-\varepsilon} (y) \, \, \delta_{\varepsilon}(y) \ = \ y$$
for any $y$ and for any $\varepsilon \in (0,1)$. This fact implies that $\displaystyle \left\{ \delta_{\mu} y \mbox{ : } \mu \in (0,+\infty) \right\}$
is a one parameter semigroup. Moreover, let $\displaystyle f_{y}: \mathbb{R} \rightarrow X$, defined by: if $\varepsilon > 0$ then 
$\displaystyle f_{y} (\varepsilon) = \delta_{\varepsilon} y$, else $\displaystyle f_{y}(\varepsilon) = \delta_{\varepsilon} y^{-1}$. Then 
$\displaystyle f_{y}$ is a group morphism from $\mathbb{R}$ to $X$, with $\displaystyle f_{y}(1) = y$. Therefore $\displaystyle f_{y} (\varepsilon) \ = \ 
\exp(\varepsilon y) \ = \ \varepsilon y$. According to definition \ref{defnormedhom} the group $X$ is identified with its Lie algebra and any element $y$ has a decomposition $\displaystyle y = y_{1} + y_{2} + ... + y_{m}$ and $\displaystyle \delta_{\varepsilon} y \ = \ \sum_{j=1}^{m} \varepsilon^{j} y_{j}$. We 
proved that $m=1$, otherwise said that the graduation of the group has only one level, that is the group is abelian. \quad $\square$

\subsection{The ratio of three collinear points}

In this section we prove that the noncommutative affine geometry is a geometry 
in the sense of the Erlangen program, because it can be described as the geometry 
of collinear triples (see definition \ref{colin}). Collinear triples generalize the 
basic ratio invariant of classical affine geometry. 

Indeed, theorem \ref{thmenelin} provides us with a mean to introduce a version of 
the ratio of three collinear points in a strong linear dilatation structure. We define here 
 {\it collinear triples},  {\it the ratio function} and {\it the ratio norm}. 

\begin{dfn}
Let $(X, d, \delta)$ be a strong linear dilatation structure.  Denote by 
$\displaystyle x^{\alpha} = (x, \alpha)$, for any $x \in X$ and $\alpha \in (0,
+ \infty)$.  An ordered set 
$\displaystyle (x^{\alpha}, y^{\beta}, z^{\gamma} )\in 
\left(X \times (0,+\infty)\right)^{3}$ is a collinear triple if: 
\begin{enumerate}
\item[(a)] $\alpha \beta \gamma = 1$ and all three numbers 
are different from $1$, 
\item[(b)] we have $\displaystyle \delta_{\alpha}^{x} \, \delta^{y}_{\beta} 
\, \delta^{z}_{\gamma} = \, id$. 
\end{enumerate}
The ratio norm $\displaystyle r(x^{\alpha},y^{\beta},z^{\gamma})$ of the collinear triple 
$\displaystyle (x^{\alpha},y^{\beta},z^{\gamma})$  is given
by the expression: 
$$r(x^{\alpha},y^{\beta},z^{\gamma}) \, = \, \frac{\alpha}{1 - \alpha \beta}$$
 Let $\displaystyle 
(x^{\alpha},y^{\beta},z^{\gamma})$ be a collinear triple. Then we have: 
$\displaystyle \delta^{x}_{\alpha} \, \delta_{\beta}^{y} \, = \, \delta^{z}_{\alpha \beta}$  
with $\alpha, \beta, \alpha \beta$ not equal to $1$. By theorem \ref{thmenelin} 
the point $z$ is uniquely determined by $\displaystyle (x^{\alpha}, y^{\beta})$,
 therefore we can express it as a function 
 $z = w(x,y, \alpha, \beta)$. 
The function $w$ is called the ratio function. 
\label{colin}
\end{dfn}

In the next proposition we obtain a formula for $w(x,y, \alpha, \beta)$. Alternatively this can be seen as 
another proof of theorem \ref{thmenelin}. 

\begin{prop}
In the hypothesis of proposition \ref{dilinv}, for any $\varepsilon, \mu  \in (0,1)$  and $x, y \in X$ 
we have: 
$$w(x,y, \varepsilon, \mu) \ = \ g_{\varepsilon \mu} \left( h_{\varepsilon}(x) h_{\mu} (\delta_{\varepsilon} y ) \right) $$
\end{prop}

\paragraph{Proof.}
Any $z \in X$ with the property that for any $u \in X$ we have $\displaystyle \delta_{\varepsilon}^{x} \, \delta_{\mu}^{y} (u) \ = \ \delta^{z}_{\varepsilon \mu} (u)$ satisfies the equation: 
\begin{equation}
x \, \delta_{\varepsilon} \left( x^{-1} y \delta_{\mu} (y^{-1}) \right) \ = \ z \delta_{\varepsilon \mu} (z^{-1})
\label{equ}
\end{equation}
This equation can be put as: 
$$h_{\varepsilon}(x) \, \delta_{\varepsilon} \left( h_{\mu}(y) \right) \ = \ h_{\varepsilon \mu} (z)$$\
From proposition \ref{dilinv} we obtain that indeed exists and it is unique $z \in X$ solution of this 
equation. We use further homogeneity 
of $\displaystyle h_{\mu}$ and we get: 
$$z \ = \ w(x,y, \varepsilon, \mu) \ = \ g_{\varepsilon \mu} \left( h_{\varepsilon}(x) h_{\mu} (\delta_{\varepsilon} y ) \right)  \quad \square $$

Remark that if $\displaystyle (x^{\alpha},y^{\beta},z^{\gamma})$ is a collinear triple then any circular permutation 
of the triple is also a collinear triple. We can not deduce from here a collinearity notion for the triple of {\it points} $\left\{ x, y, z \right\}$. 
Indeed, as the following example shows, even if $\displaystyle (x^{\alpha},y^{\beta},z^{\gamma})$ is a collinear triple, it may happen that 
here are no numbers  $\alpha', \beta', \gamma'$ such that $\displaystyle (y^{\beta'},x^{\alpha'},z^{\gamma'})$ is a collinear triple.

\paragraph{Collinear triples in the Heisenberg group.} The Heisenberg group $H(n) = \mathbb{R}^{2n+1}$ is a 2-step  Carnot group. For the points 
of $X \in H(n)$ we use the notation $X = (x, \bar{x})$, with $\displaystyle x \in \mathbb{R}^{2n}$ and $ \bar{x} \in \mathbb{R}$.   The  group 
operation is : 
$$X \, Y \ = \ (x,\bar{x})  (y,\bar{y}) \  = \ (x + y, \bar{x} + \bar{y} + \frac{1}{2} \omega(x,y))$$
where $\omega$ is the standard symplectic form on $\mathbb{R}^{2n}$. We shall identify 
the Lie algebra with the Lie group. The bracket is 
$$[(x,\bar{x}),(y,\bar{y})] = (0, \omega(x,y))$$
The Heisenberg algebra is generated by 
$$V = \mathbb{R}^{2n} \times \left\{ 0 \right\}$$ 
and we have the relations $V + [V,V] = H(n)$, $\left\{0\right\} \times \mathbb{R} \ = \ [V,V] \ = \ Z(H(n))$.

The dilatations on $H(n)$ are 
$$\delta_{\varepsilon} (x,\bar{x}) = (\varepsilon x , \varepsilon^{2} \bar{x})$$
For $X = (x, \bar{x}), Y = (y, \bar{y}) \in H(n)$ and $\varepsilon, \mu \in (0, + \infty)$, $\varepsilon \mu \not = 1$, we compute 
$Z  =  (z , \bar{z}) =  w(\tilde{x}, \tilde{y}, \varepsilon, \mu)$ with the help of equation (\ref{equ}). This equation writes: 
$$((1 - \varepsilon) x , (1 -\varepsilon^{2}) \bar{x}) \,  (\varepsilon (1- \mu) y , \varepsilon^{2}(1 -\mu^{2}) \bar{y}) \ = \ ((1- \varepsilon \mu) z, (1 - \varepsilon^{2}\mu^{2}) \bar{z})$$
After using the expression of the group operation we obtain: 
$$Z = \left( \frac{1 - \varepsilon}{1- \varepsilon \mu}  x  + \frac{\varepsilon(1 - \mu)}{1- \varepsilon \mu} y, 
 \frac{1 - \varepsilon^{2}}{1- \varepsilon^{2} \mu^{2}}  \bar{x}  + \frac{\varepsilon^{2}(1 - \mu^{2})}{1- \varepsilon^{2} \mu^{2}}  \bar{y}   + 
 \frac{\varepsilon (1-\varepsilon)(1-\mu)}{2(1- \varepsilon^{2} \mu^{2})}  \omega(x,y) \right)$$
Suppose now that $\displaystyle (X^{\alpha},Y^{\beta},Z^{\gamma})$ and $(Y^{\beta'},X^{\alpha'},Z^{\gamma'})$ are collinear triples such that 
$X = (x, 0)$, $Y = (y, 0)$ and $\omega(x,y) \not = 0$.  From the computation of the ratio function, we get that there exist numbers $k, k' \not = 0 $ such that: 
$$z \, = \, k  x \, + \, (1-k) y \, = \, (1-k') x \, + \, k' y \quad , $$ 
$$ \bar{z} \ = \ \frac{k(1-k)}{2} \, \omega(x,y) \, = \, \frac{k'(1-k')}{2} \, \omega(y,x)$$
From the equalities concerning $z$ we get that   $k'= 1-k$. This lead us to contradiction in the equalities concerning $\bar{z}$. Therefore, in this case, 
if $\displaystyle (X^{\alpha},Y^{\beta},Z^{\gamma})$ is a collinear triple then there are no $\alpha', \beta', \gamma'$ such that $(Y^{\beta'},X^{\alpha'},Z^{\gamma'})$ is a collinear triple. \quad $\square$

In a general linear dilatation structure the relation  (\ref{relation}) does not
hold. Nevertheless, there is some content of this relation which survives in 
the general context.

\begin{prop}
For $x, y$ sufficiently closed and for $\varepsilon, \mu \in \Gamma$ with $\nu(\varepsilon), \nu(\mu) \in (0,1)$,  we have the distance estimates: 
\begin{equation}
d(x, w(x,y,\varepsilon, \mu)) \ \leq \ \frac{\nu(\varepsilon)}{1 - \nu( \varepsilon \mu)} d( x, \delta^{y}_{\mu} x)   
\label{equa1}
\end{equation}
\begin{equation}
d(y, w(x,y,\varepsilon, \mu)) \ \leq \ \frac{1}{1 - \nu(\varepsilon \mu)} d( y, \delta^{x}_{\varepsilon} y)   
\label{equa2}
\end{equation} 
\label{plimit}
\end{prop}

\paragraph{Proof.} 
Further we shall use the notations from the proof of theorem \ref{thmenelin}, in particular $w = w(x, y, \varepsilon, \mu)$. We define by induction four sequences 
of points (the first two sequences are defined as in relation (\ref{definduc})): 
$$x_{n+1} \ = \ \delta_{\mu}^{\delta_{\varepsilon}^{x_{n}} y_{n}} x_{n} \quad , 
\quad y_{n+1} \ = \ \delta_{\varepsilon}^{x_{n}} y_{n} $$
$$x_{n+1}' \ = \ \delta_{\mu}^{\delta_{\varepsilon}^{y_{n}'} x_{n}'} x_{n} \quad , 
\quad y_{n+1}' \ = \ \delta_{\varepsilon}^{x_{n+1}'} y_{n}' $$
with initial conditions $\displaystyle x_{0} = x, y_{0} = y, x_{0}' = x, y_{0}' = \delta_{\varepsilon}^{x} y$. 
The first two sequences are like in the proof of theorem \ref{thmenelin}, while for the third and fourth sequences we have the relations 
$\displaystyle x_{n}' = x_{n}$, $\displaystyle y_{n}' = y_{n+1}$. These last sequences  come from the fact that they appear if we repeat the proof 
of theorem \ref{thmenelin} starting from the relation: 
$$\delta^{\delta^{x}_{\varepsilon} y}_{\mu} \delta^{x}_{\varepsilon} \ = \ \delta^{w}_{\varepsilon \mu}$$

  We know that all these four sequences converge to $w$ as $n$ goes to $\infty$. 
Moreover, we know from the proof of theorem \ref{thmenelin} that for all $n \in \mathbb{N}$ we have 
$$\displaystyle d(x_{n}, x_{n+1}) \ = \ d(x, \delta^{x}_{\varepsilon} \delta^{y}_{\mu} x) \nu(\varepsilon \mu )^{n}$$
There is an equivalent relation in terms of the sequence $\displaystyle y_{n}'$, which is the following:
$$\displaystyle d(y_{n}', y_{n+1}') \ = \ d(\delta^{x}_{\varepsilon} y, \delta^{\delta^{x}_{\varepsilon} y}_{\mu} \delta^{x}_{\varepsilon} \delta^{x}_{\varepsilon} y) \nu(\varepsilon \mu )^{n}$$
This relation becomes: for any $n \in \mathbb{N}$, $n \geq 1$ 
$$d(y_{n}, y_{n+1}) \ = \ d(y, \delta^{x}_{\varepsilon} y) \nu(\varepsilon \mu)^{n+1}$$ 
For the first distance estimate we write: 
$$d(x, w) \ \leq \ \sum_{n=0}^{\infty} d(x_{n}, x_{n+1}) \ = \  d(x, \delta^{x}_{\varepsilon} \delta^{y}_{\mu} x) \, \left(\sum_{n = 0}^{\infty} \nu(\varepsilon \mu)^{n} \right) \ = \ \frac{\nu(\varepsilon)}{1 - \nu(\varepsilon \mu)} d(x,\delta^{y}_{\mu} x) $$
For the second distance estimate we write: 
$$d(y, w) \ \leq \ d(y, y_{1}) + \sum_{n=1}^{\infty} d(y_{n}, y_{n+1}) \ = \ 
d(y, y_{1}) + \frac{\nu(\varepsilon \mu)}{1 - \nu(\varepsilon \mu)} d(y, \delta^{x}_{\varepsilon} y) \ = $$
$$= \  d(y , \delta^{x}_{\varepsilon}y) \left( 1 + \frac{\nu(\varepsilon \mu)}{1 - \nu(\varepsilon \mu)} \right) \ =  \ 
\frac{1}{1 - \nu(\varepsilon \mu)} d( y, \delta^{x}_{\varepsilon} y) $$
and the proof is done. \quad $\square$

For a collinear triple $\displaystyle (x^{\alpha},y^{\beta},z^{\gamma})$ in a 
general linear dilatation structure we cannot say that $x, y, z$ lie on
the same geodesic. This is  false, as shown by easy examples in the 
Heisenberg group, the simplest noncommutative Carnot group. 

Nevertheless, theorem \ref{thmenelin} allows to speak about collinearity in the 
sense of definition \ref{colin}.

Affine geometry is the study of relations which are invariant with respect to the
group of affine transformations. An invertible transformation is affine 
if and only if it preserves the ratio of any three collinear points. We are 
thus arriving  to the following definition. 

\begin{dfn}
Let $(X,d,\delta)$ be a linear dilatation structure. A geometrically  affine transformation 
 $T: X \rightarrow X$ is a Lipschitz invertible transformation such that 
 for any collinear triple $\displaystyle (x^{\alpha},y^{\beta},z^{\gamma})$ the 
 triple $\displaystyle ((Tx)^{\alpha},(Ty)^{\beta},(Tz)^{\gamma})$ is collinear. 
\label{afftran}
\end{dfn}

The group of geometric affine transformations defines a geometry in the sense of Erlangen program. 
The main invariants of such a geometry are collinear triples. There is no obvious connection between
collinearity and geodesics of the space, as in classical affine geometry. (It is 
worthy to notice that in fact, there might be no geodesics in the metric space 
$(X,d)$ of the linear dilatation structure $(X,d,\delta)$. For example, there
are linear dilatation structures defined over the boundary of the dyadic tree 
\cite{buligaself}, which is homeomorphic with the middle thirds Cantor set.)

The first result for such a geometry is the following. 

\begin{thm}
Let $(X,d,\delta)$ be a strong linear dilatation structure. Any Lipschitz, invertible, 
transformation $T: (X,d) \rightarrow (X,d)$ is affine in the sense of  
definition \ref{defgl} if and only if it is geometrically affine in the sense of 
definition \ref{afftran}. 
\label{tafflin}
\end{thm}

\paragraph{Proof.} 
The first implication, namely $T$ affine in the sense of  
definition \ref{defgl} implies $T$ affine in the sense of 
definition \ref{afftran}, is 
straightforward: by hypothesis on $T$, for any collinear 
  triple  $\displaystyle (x^{\alpha},y^{\beta},z^{\gamma})$ we have  the
  relation 
$$\displaystyle T \, \delta_{\alpha}^{x} \, \delta^{y}_{\beta} 
\, \delta^{z}_{\gamma} \, T^{-1} \, = \, \delta_{\alpha}^{Tx} \, 
\delta^{Ty}_{\beta} 
\, \delta^{Tz}_{\gamma}$$  
Therefore, if $\displaystyle (x^{\alpha},y^{\beta},z^{\gamma})$ is a collinear 
triple then the triple  $\displaystyle ((Tx)^{\alpha},(Ty)^{\beta},(Tz)^{\gamma})$ is 
collinear. 

In order to show the inverse implication we use the linearity of the dilatation structure. Let $x, y \in X$ and 
$\varepsilon, \eta \in \Gamma$. Then 
$$\delta^{x}_{\varepsilon} \, \delta^{y}_{\eta} \delta^{x}_{\varepsilon^{-1}} \ = \ \delta^{\delta^{x}_{\varepsilon} y}_{\eta}$$
This identity leads us to the description of $\displaystyle \delta^{x}_{\varepsilon} y$ in terms of the ratio function. Indeed, we have: 
$$\delta_{\varepsilon}^{x} y \ = \ w(w(x, y, \varepsilon, \eta), \varepsilon \eta, \varepsilon^{-1})$$
If the transformation $T$ is geometrically affine then we easily find that it is affine in the sense of definition \ref{defgl}: 
$$T \left( \delta_{\varepsilon}^{x} y \right) \ = \ w(w(T x, T y, \varepsilon, \eta), \varepsilon \eta, \varepsilon^{-1}) \ = \  \delta^{Tx}_{\varepsilon} Ty$$ 
 \quad $\square$

As a  conclusion for this section, theorem \ref{tafflin} shows that  in a linear
dilatation structure we may take dilatations as the basic affine invariants. It 
is surprising that in such a geometry there is no obvious notion of a line, due 
to the fact that not any permutation of a collinear triple is again a collinear
triple.

\end{document}